\documentclass[11pt]{article}

\usepackage{color}
\usepackage{latexsym}
\usepackage{amssymb}
\usepackage{amsthm}
\usepackage{amsmath}
\usepackage{graphicx}
\usepackage{enumerate}
\usepackage{hyperref}
\usepackage{color}

\newtheorem{Theorem}{Theorem}[part]
\newtheorem{Definition}{Definition}[part]
\newtheorem{Proposition}{Proposition}[part]

\newtheorem{Remark}{Remark}[part]

\def \N{\mathbb{N}}
\def \R{\mathbb{R}}

\def \E{\mathbb{E}}
\def \F{\mathbb{F}}
\def \G{\mathbb{G}}

\def \P{\mathbf{P}}

\def \Bc{{\cal B}}
\def \Cc{{\cal C}}

\def \Fc{{\cal F}}
\def \Gc{{\cal G}}

\def \Ic{{\cal I}}

\def \Pc{{\cal P}}

\def \Sc{{\cal S}}
\def \Tc{{\cal T}}

\def \Eb{{\bf E}}
\def \Pb{{\bf P}}

\def \ni{\noindent}

\def \eps{\varepsilon}

\def \ep{\hbox{ }\hfill$\Box$}

\def\reff#1{{\rm(\ref{#1})}}

\def\beqs{\begin{eqnarray*}}
\def\enqs{\end{eqnarray*}}
\def\beq{\begin{eqnarray}}
\def\enq{\end{eqnarray}}

\definecolor{darkred}{rgb}{0.8,0,0}
\definecolor{darkblue}{rgb}{0,0,0.7}
\definecolor{darkgreen}{rgb}{0,0.4,0}

\begin{document}
\title{Adding constraints to BSDEs with Jumps: \\ an alternative to multidimensional reflections
\vspace{3mm}}
\date{ January 2011}
\author{
\begin{tabular}{ccc}
Romuald ELIE $\&$ Idris KHARROUBI\\
  \small{CEREMADE, CNRS, UMR 7534, } \\
 \small{Universit{\'e} Paris-Dauphine,}\\
\small{and CREST} \\
\small{\sf $\{$elie,kharroubi$\}$@ceremade.dauphine.fr }
\end{tabular}
\vspace{5mm}
}

\maketitle

\begin{center}
\vspace{-9mm}
\end{center}
 

\begin{abstract} 

This paper is dedicated to the analysis of backward stochastic differential equations (BSDEs) with jumps, subject to an additional global constraint involving all the components of the solution. We study the existence and uniqueness of a minimal solution for these so-called constrained BSDEs with jumps via a penalization procedure. This new type of BSDE offers a nice and practical unifying framework 
to the notions of  constrained BSDEs presented in  \cite{penxu07} and BSDEs with constrained jumps introduced in \cite{kmpz08}. More remarkably, the solution of a multidimensional Brownian reflected BSDE studied in  \cite{hz08} and \cite{ht07} can also be represented via a well chosen one-dimensional constrained BSDE with jumps.
This last result is very promising from a numerical point of view for the resolution of high dimensional optimal switching problems and more generally for systems of coupled variational inequalities.

\end{abstract}

\vspace{5mm}


\vspace{5mm}

\noindent {\bf MSC Classification (2000):}  93E20, 60H30, 60J75.




\section{Introduction}
\setcounter{equation}{0} \setcounter{Assumption}{0}
\setcounter{Theorem}{0} \setcounter{Proposition}{0}
\setcounter{Corollary}{0} \setcounter{Lemma}{0}
\setcounter{Definition}{0} \setcounter{Remark}{0}

Since their introduction by Pardoux and Peng in \cite{parpen90}, Backward Stochastic Differential Equations (BSDEs in short) have been widely studied. In particular, they appear as a very powerful tool to solve partial differential equations (PDEs) and corresponding stochastic optimization problems. Several generalizations of this notion are based on the addition of new constraints on the solution. First, El Karoui et al. \cite{elketal97} study the case where the component $Y$ is forced to stay above a given process, leading to the notion of reflected BSDEs related to optimal stopping and obstacle problems. Motivated by super replication issues under portfolio constraints, Cvitanic et al. \cite{cvikarson98}  consider the case where the component $Z$ is constrained to stay in a fixed convex set. More recently, Kharroubi et al. \cite{kmpz08} introduce a constraint on the jump component $U$ of the BSDE, providing a representation of solutions for a class of PDE, called quasi-variational inequalities, arising from optimal impulse control problems. The generalization of the results of  El Karoui et al. \cite{elketal97} to oblique reflections in a multi-dimensional framework was first given in a very special case (e.g. the generator does not depend on z) by 
Ramasubramanian \cite{ram02}, who studied a BSDE reflected in an orthant. Then, 
Hu and Tang \cite{ht07} followed by Hamad\`ene and Zhang \cite{hz08} consider general BSDEs with oblique reflections and connect them with systems of variational inequalities and optimal switching problems. 
Our paper introduces the notion of constrained BSDEs with jumps, which offers in particular a nice and natural probabilistic representation for these types of switching problems. This new notion essentially unifies and extends the notions of constrained BSDE without jumps, BSDE with constrained jumps as well as multidimensional BSDE with oblique reflections.\\

Let us illustrate our presentation with the example of the following switching problem 
\beq \label{switch}
 \sup_{\alpha} \Eb\Big[ g_{\alpha_T}(X_T) + \int_0^T \psi_{\alpha_s}(s,X_s) ds + \sum_{0 < \tau_k\leq T} c_{\alpha_{\tau_k^-},\alpha_{\tau_k}} \Big]\,,
 \enq
 where $X$ is an underlying It\^o diffusion process, $\alpha$ is a switching control process valued in $\Ic:=\{1,\ldots,m\}$, $m>0$, and $(\tau_k)_k$ denotes the jump times of the control $\alpha$. 
This type of stochastic control problem is typically encountered by an agent maximizing the production rentability of a given good  by switching between $m$ possible modes of production based on different commodities. A switch is penalized by a given cost function $c$ and the production rentability functions $\psi$ and $g$ depend on the  chosen mode of production.   
 As observed in \cite{djehampop07}, the solution of problem \reff{switch} starting in mode $i_0\in\Ic$ at time $t$ rewrites $Y^{i_0}_t$ where $(Y^i,Z^i,K^i)_{i\in\Ic}$ solves the following multidimensional reflected BSDE
 \begin{equation}\label{BSDEobliqueintro}
\left\{
\begin{array}{l} \vspace{1mm}
Y^i_t  =  g_i(X_T)  + \int_t^T \psi_i(s,X_s)ds - \int_t^T \langle Z^i_s,dW_s\rangle + K^i_T - K^i_t \,, \quad 0\le t \le T\,, \quad i\in\Ic \\ \vspace{1mm}
Y^i_t  \geq  Y^j_t + c_{i,j}   \,,  \quad 0\le t \le T\,,  \quad i,j \in\Ic \;,   \\ 
\int_0^T[Y_t^i-\max_{j\in \Ic} \{ Y_j + c_{i,j} \}] dK^i_t=0 \;, \quad i\in\Ic \;.
\end{array}
\right.
\end{equation}
The main difficulty in the derivation of a one-dimensional BSDE representation for this type of problem relies on the dependence of the solution in mode $i\in\Ic$ with respect to the global solution in all possible modes. Nevertheless, Tang and Yong \cite{ty93} interpret the value function associated to this problem  as the unique viscosity solution of a given coupled system of variational inequalities.  A clever observation of Bouchard \cite{bou06} concludes that this unique viscosity solution represents also the value function of a well suited stochastic target problem associated to a diffusion with jumps. Using entirely probabilistic arguments, the BSDE representation provided in this paper heavily relies on this type of correspondence. In our approach, we let artificially the strategy jump randomly between the different modes of production. Similarly to the approach of Pardoux et al. \cite{ppr97}, this allows to retrieve in the jump component of a one-dimensional backward process, some information regarding the solution in the other modes of production. Indeed, let us introduce a pure jump process $(I_t)_{0\leq t\leq T}$ based on an independent random measure $\mu$ and 
consider the following constrained BSDE associated to the two dimensional forward process $(I,X)$ (called transmutation-diffusion process in \cite{ppr97}) and defined on $[0,T]$ by:
\begin{equation}\label{BSDEgenintro} 
\left\{
 \begin{array}{lcl}
  \tilde Y_t \!&=\!\!& g_{I_T}(X_T) + \int_t^T \psi_{I_s}(s,X_s) ds  + \tilde K_T - \tilde K_t  -  \int_t^T \langle \tilde Z_s, dW_s\rangle  -  \int_t^T\!\int_\Ic  \tilde U_s(i)   \mu(ds,di)\,,\\
  \tilde U_t(i) \!\! \!& \geq \!\!& c_{i,I_{t-}}, \;\;\;\;\; d\Pb\otimes dt\otimes\lambda(di) \;\;  a.e.
 \end{array}
 \right.
\end{equation}
This BSDE enters into the class of constrained BSDEs studied in the paper and its unique minimal solution relates directly to the solution of \reff{BSDEobliqueintro} via the relation  $( \tilde Y_t, \tilde Z_t, \tilde U_t) = (Y^{I_t}_t, Z^{I_t}_t, \{Y^i_t-Y^{I_{t-}}_{t^-}\}_{i\in\Ic})$ for $t\in[0,T]$.
 In particular, the solution of the switching problem \reff{switch} starting in mode $I_0$ at time $0$ rewrites $\tilde Y^{I_0}_0$.\\

 In order to unify our results with the one based on multidimensional reflected BSDE considered in \cite{ht07} or \cite{hz08}, we extend this approach and introduce the notion of constrained BSDE with jumps whose solution $(Y,Z,U,K)$ satisfies the general dynamics
  \beq
 Y_t &=& \xi + \int_t^T f(s,Y_s,Z_s,U_{s}) ds  + K_T - K_t  -  \int_t^T \langle Z_s, dW_s\rangle   \label{BSDEgenIntro}  
 -  \int_t^T\int_\Ic  U_s(i)   \mu(ds,di),  \qquad
 \enq
 a.s., for $0 \leq t \leq T$,  as well as the constraint
 \beq \label{hconsIntro}
 h_i(t,Y_{t^-},Z_{t},U_t(i)) & \geq & 0, \;\;\;\;\; d\Pb\otimes dt\otimes\lambda(di) \;\;  a.e.\; ,
 \enq
 where $f$ and $h$ are given random Lipschitz functions, and $h$ is non-increasing in its last variable. Through a penalization argument, we provide in Section 2 the existence of a unique minimal solution to  the constrained BSDE with jumps \reff{BSDEgenIntro}-\reff{hconsIntro}. This new type of BSDE 
 mainly extends and unifies the existing literature on BSDEs in three interconnected directions:
 \begin{itemize}
 \item We generalize the notion of BSDE with constrained jumps considered in \cite{kmpz08}, letting the driver function $f$ depend on $U$ and considering a general constraint function $h$ depending on all the components of the solution.
 \item We add some jumps in the dynamics of constrained BSDE studied in \cite{penxu07} and let the coefficients depend on the jump component $U$.
 \item Via the addition of artificial jumps, a well chosen one-dimensional constrained BSDE with jumps allows to represent the solution of a multidimensional reflected BSDE, in the framework of \cite{hz08} or \cite{ht07}.
 \end{itemize}
 
 We believe that the representation of a multidimensional obliquely reflected BSDE by a one-dimensional constrained BSDE with jumps is {also} numerically very promising. 
 As developed in \cite{bptw11}, it offers the possibility to solve high-dimensional optimal switching problems via a natural extension of the entirely probabilistic numerical scheme studied in  \cite{be06}. 
Such type of algorithm could also solve high dimensional systems of variational inequalities, which relate directly to multidimensional BSDEs with oblique reflections, see \cite{ht07} for more details. The algorithm as well as the Feynman Kac representation of general constrained BSDEs with jumps are presented in \cite{ek09}. \\

 The paper is organized as follows. The next section provides the existence of a unique minimal solution for the new class of constrained BSDEs with jumps  \reff{BSDEgenIntro}-\reff{hconsIntro}. The connection  with multidimensional reflected BSDEs is detailed in Section 3.  We regroup in the last section of the paper some technical results on BSDEs, mainly extensions of existing results, which are not the main focus of the paper but present some interest in themselves: we provide a comparison and a monotonic limit theorem for reflected BSDEs with jumps, as well as viability and comparison properties for multidimensional constrained BSDEs. We isolate these results in order to present them in a general framework and to simplify their possible future invocation.
All the proofs of the paper only rely on probabilistic arguments and can be applied in a non-Markovian setting.

\paragraph{Notations.} Throughout this paper we are given a finite terminal time $T$
and a probability space $(\Omega,\Gc,\Pb)$ endowed with a $d$-dimensional standard Brownian
motion $W$ $=$ $(W_t)_{t\geq 0}$, and a Poisson random
measure $\mu$ on $\R_+\times \Ic$, where $\Ic$ $=$
$\{1,\ldots,m\}$,  with intensity measure $\lambda(di)dt$ for some
finite measure $\lambda$ on $\Ic$ with $\lambda(i)$ $>$ $0$ for all $i$ $\in$ $\Ic$. We set $\tilde\mu(dt,di)$ $=$
$\mu(dt,di)-\lambda(di)dt$ the compensated measure associated to
$\mu$. $\sigma(\Ic)$ denotes the $\sigma$-algebra of subsets of $\Ic$. For $x$ $=$ $(x_{1},\ldots,x_{\ell})$ $\in$ $\R^\ell$ with $\ell$ $\in$ $\N$, we set $|x|$ $=$ $\sqrt{|x_{1}|^2+\cdots+|x_{\ell}|^2}$ the Euclidean norm. We denote by $\G$ $=$ $(\Gc_t)_{t\geq 0}$ (resp. $\F$ $=$ $(\Fc_t)_{t\geq 0}$) the augmentation of
the natural filtration generated by $W$ and $\mu$ (resp. by $W$),  and by $\Pc_{\G}$ (resp. $\Pc_{\F}$, $\mathfrak{P}_{\G}$, $\mathfrak{P}_{\F}$) the
$\sigma$-algebra of $\G$-predictable (resp. $\F$-predictable $\G$-progressive, $\F$-progressive) subsets of $\Omega\times [0,T]$.
 We denote
 by ${\bf \Sc^2_{\mathbb{G}}}$ (resp.  ${\bf \Sc^{2}_{\mathbb{F}}}$) 
 the set of real-valued  c\`ad-l\`ag $\mathbb{G}$-adapted (resp. continuous $\mathbb{F}$-adapted) 
processes $Y$ $=$ $(Y_t)_{0\leq t\leq T}$ such that
 \beqs
\|Y\|_{_{{\bf  \Sc^2}}} \;:=\; \left(\Eb\Big[ \sup_{0\leq t\leq T} |Y_t|^2 \Big]\right)^{1\over 2}  \; < \; \infty.
 \enqs
${\bf L^p(0,T)}$, $p$ $\geq$ $1$,  is the set of real-valued processes $\phi$ $=$ $(\phi_t)_{0\leq t\leq T}$ such that
 \beqs
\|\phi\|_{_{\bf L^p(0,T)}} ~ := ~ \Big(\Eb\Big[\int_0^T |\phi_t|^p dt\Big] \Big)^{1 \over p}\;<\; \infty,
 \enqs
and ${\bf L^p_{\F}(0,T)}$ (resp. ${\bf L^p_{\G}(0,T)}$) is the subset of ${\bf L^p(0,T)}$ consisting of $\mathfrak{P}_{\F}$-measurable  (resp. $\mathfrak{P}_{\G}$-measurable) processes.
 ${\bf L^p_{\F}(W)}$ (resp. ${\bf L^p_{\G}(W)}$), $p$ $\geq$ $1$, is the set of
$\R^d$-valued $\Pc_{\F}$-measurable  (resp. $\Pc_{\G}$-measurable) processes $Z$ $=$ $(Z_t)_{0\leq t\leq
T}\in{\bf L^p_{\F}(0,T)}$ (resp. ${\bf L^p_{\G}(0,T)}$)  .
${\bf L^p(\tilde\mu)}$, $p$ $\geq$ $1$,  is the set of
$\Pc\otimes\sigma(\Ic)$-measurable maps $U$ $:$ $\Omega\times [0,T]\times
\Ic$ $\rightarrow$ $\R$ such that \beqs \| U\|_{_{{\bf
L^p(\tilde\mu)}}} \;:= \;  \left(\Eb\Big[ \int_0^T\int_\Ic  |U_t(i)|^p
\lambda(di) dt \Big]\right)^{1\over p}  \; < \; \infty.
\enqs
${\bf A^2_{\F}}$ (resp. ${\bf A^2_{\G}}$) is the  closed subset of ${\bf \Sc^2_{\F}}$ (resp. ${\bf \Sc^2_{\G}}$) consisting of
nondecreasing  processes $K$ $=$ $(K_t)_{0\leq t\leq T}$ with $K_0$
$=$ $0$.
Finally, for $t\in[ 0,T]$, $\Tc_t$ denotes the set of $\F$-stopping times $\tau$ such that  $\tau \in[t,T]$, $\P$-a.s.. For ease of notation, we omit in all the paper the dependence in $\omega\in\Omega$,  whenever it is not relevant.


\section{Constrained Backward SDEs with jumps}
\setcounter{equation}{0} \setcounter{Assumption}{0}
\setcounter{Theorem}{0} \setcounter{Proposition}{0}
\setcounter{Corollary}{0} \setcounter{Lemma}{0}
\setcounter{Definition}{0} \setcounter{Remark}{0}

This section is devoted to the presentation of constrained Backward SDEs with jumps, generalizing the framework considered in \cite{kmpz08} or \cite{penxu07}. Namely: 
\begin{itemize}
 \item We allow the driver function to depend on the jump component of the backward process,  
 \item We extend the class of possible constraint functions by letting them depend on all the components of the solution to the BSDE. 
 \end{itemize}

 We adapt the arguments developed in \cite{kmpz08} in order to derive existence and uniqueness of a minimal solution for this new type of BSDE. No major difficulty appears for the obtention of these results and, from our point of view, the nice feature of such constrained BSDE relies on their relation with multidimensional reflected BSDE, developed in the next section.  In order to simplify, the readability of the paper, the required technical extensions of comparison and monotonic limit theorems are reported in Sections \ref{App Comp Jumps} and \ref{App Monotonic Jumps}. They are presented in a more abstract framework and can therefore be quoted more conveniently in the future.


\subsection{Formulation}

A constrained BSDE with jumps is characterized by three objects:

\begin{itemize}
\item a terminal condition, i.e. a $\Gc_{T}$-measurable random variable $\xi$,

\item  a driver function, i.e. a map $f$ $:$ $\Omega\times[0,T]\times\R\times\R^d\times\R^m \rightarrow\R$,  which is $\mathfrak{P}_{\G}\otimes \Bc(\R)\otimes \Bc(\R^d)\otimes \Bc(\R^m)$-measurable,

\item a constraint function, i.e. a $\sigma(\Ic)\otimes\mathfrak{P}_{\G}\otimes\Bc(\R)\otimes\Bc(\R^d)\otimes\Bc(\R)$-measurable map $h~:~\Ic\times\Omega\times[0,T]\times\R\times\R^d\times\R$ $\rightarrow$ $\R$  such that $h_i(\omega,t,y,z,.)$ is non-increasing for all $(i,\omega,t,y,z)$ $\in$ $\Ic\times\Omega\times[0,T]\times\R\times\R^d$.
\end{itemize}

\begin{Definition}
(i) A solution to the corresponding constrained BSDE with jumps is a quadruple $(Y,Z,U,K)$ $\in$ ${\bf \Sc_{\G}^2}\times{\bf L^2_{\G}(W)}\times{\bf L^2(\tilde \mu)}\times{\bf A^2_{\G}}$ satisfying
 \beq
 Y_t &=& \xi + \int_t^T\!\!\! f(s,Y_s,Z_s,U_{s}) ds  + K_T - K_t  -  \int_t^T \!\!\langle Z_s, dW_s\rangle   \label{BSDEgen}  
 -  \int_t^T\!\int_\Ic  U_s(i)   \mu(ds,di), \quad  \qquad
 \enq
 for $0 \leq t \leq T$ a.s., as well as the constraint
 \beq \label{hcons}
 h_i(t,Y_{t^-},Z_{t},U_t(i)) & \geq & 0, \;\;\;\;\; d\Pb\otimes dt\otimes\lambda(di) \;\;  a.e.\; .
 \enq
(ii) $(Y,Z,U,K)$ is a minimal solution to \reff{BSDEgen}-\reff{hcons} whenever it is solution to \reff{BSDEgen}-\reff{hcons} and  for any other solution $(\check Y,\check Z,\check U,\check K)$  of
\reff{BSDEgen}-\reff{hcons}, we have $Y  \leq \check Y$  a.s. 
\end{Definition}
We notice that for a minimal solution $(Y,Z,U,K)$ to \reff{BSDEgen}-\reff{hcons}, the component $Y$ naturally interprets in the terminology of Peng \cite{pen99} as the smallest supersolution to \reff{BSDEgen}-\reff{hcons}. \\

\begin{Remark}
{\rm  In the case where the driver function $f$ does not depend on $U$ and the constraint function $h$ is of the form $h_i(u+c(t,y,z))$, observe that this BSDE exactly fits in the framework  considered in \cite{kmpz08}. 
Similarly, in the Brownian case (i.e. no jump component), this type of BSDEs was studied in \cite{penxu07}. Therefore, our framework generalizes and unifies those considered in  \cite{kmpz08} and \cite{penxu07}.
}
\end{Remark}

In order to work on this class of BSDE, we require the classical Lipschitz and linear growth conditions on the coefficients, as well as a control on the way the driver  function depends on the jump component $U$ of the BSDE. We regroup these conditions in the following assumption.\\

\noindent \textbf{(H0)}
 \begin{enumerate}[(i)]
\item  There exists a constant $k>0$ such that the functions $f$ and $h$ satisfy  $\P$-a.s. the uniform Lipschitz property:
\beqs
 |f(t,y,z,u)-f(t,y',z',u')|  & \leq & k |(y,z,u)-(y',z',u')|\;, \\
  |h_i(t,y,z,u_{i})-h_i(t,y',z',u_{i}')| & \leq & k   |(y,z,u_{i})-(y',z',u'_{i})|\;,
 \enqs 
  for all  
 $\{i,t,(y,z,u),(y',z',u')\}\in\Ic\times[0,T]\times[\R\times\R^{d}\times\R^m]^2$.

\item The coefficients $\xi$, $f$ and $h$ satisfy the following integrability condition
\beq
\mathbf{E}{ |\xi|^2 + \int_0^T \mathbf{E} |f(t,0,0,0)|^2dt + \sum_{i\in\Ic} \int_0^T \mathbf{E} |h_i(t,0,0,0)|^2 dt} & < & \infty\;.
\enq

\item There exist two constants $C_{1}$ $\geq$ $C_{2}$ $>$ $-1$ such that we can find a $\Pc_{\G}\otimes\sigma(\Ic)\otimes\Bc(\R)\otimes\Bc(\R^d)\otimes\Bc(\R^m)\otimes\Bc(\R^m)$-measurable map
$\gamma~:\Omega\times[0,T]\times\Ic\times\R\times\R^d\times\R^m\times\R^m\rightarrow [C_2,C_1]$ satisfying 
\beqs
f(t,y,z,u)-f(t,y,z,u') & \leq & \int_{\Ic}(u_{i}-u'_{i})\gamma_{t}^{y,z,u,u'}(i)\lambda(di),
\enqs
  for all  $(i,t,y,z,u,u')\in\Ic\times[0,T]\times\R\times\R^{d}\times[\R^m]^2$, $\P$-a.s..

\end{enumerate}

\begin{Remark}
{\rm  Under Assumption {\bf (H0)} (i) and (ii), existence and uniqueness of a solution $(Y,Z,U,K)$ to the BSDE \reff{BSDEgen} with $K$ $=$ $0$  follows from classical results on BSDEs with jumps, see Lemma 2.4 in \cite{tanli94}. In order to add the  $h$-constraint \reff{hcons}, one needs as usual to relax the dynamics of $Y$ by injecting the non-decreasing process $K$ in \reff{BSDEgen}. 
In mathematical finance, the purpose of this new process $K$ is to increase the  super replication price $Y$ of a contingent claim, under additional portfolio constraints. In order to find a minimal solution to the constrained BSDE \reff{BSDEgen}-\reff{hcons}, the nondecreasing property of $h$ is crucial for stating comparison principles needed in the penalization approach. 
}\end{Remark}

\begin{Remark}
{\rm Part (iii) of Assumption {\bf (H0)} constrains the dependence of the driver $f$ with respect to the jump component of the BSDE. It is inspired by \cite{roy06} and will ensure comparison results for BSDEs driven by this type of driver, as detailed in Section \ref{App Comp Jumps}.
}
\end{Remark}


\subsection{Approximation by penalization}\label{approxpen}

This paragraph focuses on the existence of a unique minimal solution for the constrained BSDE with jumps \reff{BSDEgen}-\reff{hcons}. Our approach requires the addition of an increasing component to the comparison results for BSDEs with jumps, derived by Royer \cite{roy06} as well as the extension of Peng's monotonic limit theorem \cite{pen99} to the consideration of BSDEs with jumps. We could not find these properties in the existing literature and report them respectively in Proposition \ref{PropComparisonBSDE} and Proposition \ref{MonThBSDEwJ} of Section \ref{SecAppendix}.  \\

The proof relies on a classical penalization argument and we introduce the following sequence of BSDEs with jumps
 \beq 
 Y_t^n &=& \xi + \int_t^T f(s,Y_s^n,Z_s^n,U_{s}^n) ds  + n \int_t^T \int_\Ic h_i^-(s,Y^{n}_{s},Z^{n}_{s},U^n_s(i)) \lambda(di)ds   \label{BSDEpen1}  \\
& & 
\qquad\qquad-  \int_t^T  \langle Z_s^n,dW_s\rangle  -  \int_t^T\int_\Ic   U_s^n(i) \mu(ds,di), \;\quad 0 \leq t \leq T,\quad n\in\N,   \nonumber \enq
 where  $h_i^-(.)$ $:=$ $\max(-h_i(.),0)$ is the negative part of the function $h_i$, $i\in\Ic$. 
 Under Assumption  {\bf (H0)}, the Lipschitz property of the coefficients $f$ and $h$ ensures existence and uniqueness of a solution $(Y^n,Z^n,U^n)$ $\in$ ${\bf \Sc_{\G}^2}\times{\bf L^2_{\G}(W)}\times{\bf L^2(\tilde\mu)}$ to \reff{BSDEpen1}, see Theorem 2.1 in \cite{barbucpar97}. \\

 In order to obtain the convergence of the sequence $(Y^n)_{n\in\N}$, we require: \\
\vspace{2mm}

\noindent {\bf (H1)} \hspace{1mm} There exists  $(\check
Y,\check Z,\check K,\check U)$ $\in$ ${\bf \Sc^2_{\G}}\times{\bf
L^2_{\G}(W)}\times{\bf L^2(\tilde \mu)}\times{\bf A^2_{\G}}$  solution of 
\reff{BSDEgen}-\reff{hcons}.\\
\vspace{2mm}

This assumption, which may appear restrictive, is rather classical and we present in Section \ref{SectionLinkReflBSDEs} a large class of cases where \textbf{(H1)} is satisfied. Furthermore, as detailed in Remark  \ref{RemUpperBoundY} below, \textbf{(H1)} can also be replaced by the weaker assumption: \\

\vspace{2mm}

\noindent {\bf (H1')} \hspace{1mm} There exists a constant $M$ such that $\sup_{n\in\N}\|Y^n\|_{\Sc^2} \le M$. \\

Under these assumptions, we are now ready to study the convergence of the quadruple $(Y^n, Z^n,U^n,K^n)_{n\in\N}$, where the nondecreasing process $K^n$ $\in$ ${\bf A^2_{\G}}$ is defined by 
 \beqs K_t^n &:=&  n \int_0^t \int_\Ic h_i^-(s,Y^n_s,Z^n_s,U^n_s(i))\lambda(di)ds, \;\;\;  0 \leq t \leq T  \;, \quad n\in\N. \enqs
 The next theorem states that the sequence $(Y^n, Z^n,U^n,K^n)_{n\in\N}$ converges indeed to the minimal solution of the constrained BSDE \reff{BSDEgen}-\reff{hcons}.

\begin{Theorem} \label{thmmain1}
Under  {\bf (H0)} and {\bf (H1)}, the following holds.
\begin{enumerate}[(i)]
\item There exists a unique minimal solution $(Y,Z,U,K)$ $\in$ ${\bf \Sc^2_{\G}}\times{\bf L^2_{\G}(W)}\times{\bf L^2(\tilde \mu)}\times{\bf A^2_{\G}}$  to \reff{BSDEgen}-\reff{hcons}, with $K$ predictable. 
\item The sequence $(Y^n)_{n\in\N}$ converges increasingly to the process $Y$ and we have
 \beqs 
 \|Y^n-Y\|_{_{{\bf L^2(0,T)}}} + 
\|Z^n-Z\|_{_{{\bf L^p(0,T)}}} + \|U^n-U\|_{_{{\bf L^p(\tilde\mu)}}} 
   & \longrightarrow_{n \rightarrow\infty} & 0,  \qquad  1\le p < 2 \;.
 \enqs 
 Moreover, $(Z,U,K)$ is the weak limit of $(Z^n,U^n,K^n)_{n\in\N}$ in ${\bf L^2_{\G}(W)}\times{\bf L^2(\tilde\mu)}\times{\bf L^2_{\G}(0,T)}$ 
and $K_{t}$ is the weak limit of $(K^n_{t})_{n\in\N}$ in ${\bf L^2(\Omega,\Gc_{t},\P)}$, for all $t\in[0,T]$. 
\end{enumerate}
\end{Theorem}

$\;$\\[-8mm]

\ni\textbf{Proof.} We prove the statements of the theorem in a reverse order. First, we show the convergence of the sequence  $(Y^n,Z^n,U^n,K^n)_{n\in\N}$. Second, we verify that the limit is a minimal solution to \reff{BSDEgen}-\reff{hcons}. Third, we tackle the uniqueness property.\\

\noindent {\bf Step 1: } {\it Convergence of $(Y^n,Z^n,U^n,K^n)_{n\in\N}$.}\\
 For $n\in\N$, we introduce the Lipschitz map $f^n$ $:=$ $f+n\int_{\Ic}h^- d\lambda$. Since $f$ satisfies {\bf (H0)}(iii) and $h$ is lipschitz and non-increasing, we deduce:
\beqs
f^n(t,y,z,u) - f^n(t,y,z,u') &\!\!\!\!\le&\!\!\!\!\! \int_{\Ic} \!\{ (u_{i}-u'_{i})\gamma_{t}^{y,z,u,u'}(i) \! +\! n(h_i^-(t,y,z,u_i) - h_i^-(t,y,z,u_i')  )   \}\lambda(di), \\
&\!\!\!\!\le&\!\!\!\!\! \int_{\Ic} (u_{i}-u'_{i}) (\gamma_{t}^{y,z,u,u'}(i) + k n {\bf 1}_{ u_i \geq u_i'}) \lambda(di), \qquad \P\mbox{- a.s} \;, \quad n\in\N\;,
\enqs
for any $(t,y,z,u,u')\in[0,T]\times\R\times\R^d\times\R^m\times\R^m$. Thus, for any $n\in\N$, the coefficients $f^n$ and $f^{n+1}$ satisfy {\bf (H0)} as well as $f^n\le f^{n+1}$. We deduce from a simplified version of Proposition \ref{PropComparisonBSDE} without the additional increasing process $K$, that the sequence $(Y^n)_{n\in\N}$ is non-decreasing.\\

Furthermore, for any quadruple $(\check Y,\check Z,\check U,\check K)$ $\in$ ${\bf \Sc^2_{\G}}\times{\bf L^2_{\G}(W)}\times{\bf L^2(\tilde \mu)}\times{\bf A^2_{\G}}$ satisfying \reff{BSDEgen}-\reff{hcons}, we obtain  $Y^n$ $\leq$ $\check Y$ a.s., $n$ $\in$ $\N$, applying once again Proposition  \ref{PropComparisonBSDE} but with coefficients $f_{1}$ $=$ $f_{2}$ $=$ $f^n$ and $K^2$ $=$ $\check K$.  
Therefore, under {\bf (H1)}, the sequence $(Y^n)_{n\in\N}$ is nondecreasing and upper bounded, ensuring its monotonic convergence to a process $Y$ with $\|Y\|_{\Sc^2}<\infty$.\\

Finally, we observe that under \textbf{(H0)} and \textbf{(H1)}, \textbf{(H3)} is satisfied by the generator $f$ and the sequence $(Y^n,Z^n,U^n,K^n)$ (see Section \ref{App Monotonic Jumps}). We are now in position to appeal to Proposition \ref{MonThBSDEwJ}, which is an extended version of Peng's monotonic limit theorem.  Hence,  the sequence $(Y^n,Z^n,U^n,K^n)_{n\in\N}$ converges in the sense specified above. Furthermore, the limit $(Y,Z,U,K)$ $\in$ ${\bf \Sc^2_{\G}}\times{\bf L^2_{\G}(W)}\times{\bf L^2(\tilde \mu)}\times{\bf A^2_{\G}}$ satisfies \reff{BSDEgen} and $K$ is predictable.\\[3mm]

\noindent {\bf Step 2: } {\it $(Y,Z,U,K)$ is a minimal solution to  \reff{BSDEgen}-\reff{hcons}.}\\
Since $(Y,Z,U,K)$ solves \reff{BSDEgen}, we now focus on the constraint property \reff{hcons}.  From the previous convergence result, we derive in particular that $(Y^n,Z^n,U^n)_{n\in\N}$ converges in ${\bf L^1_{\G}(0,T)}\times{\bf L^1_{\G}(0,T)}\times{\bf L^1(\tilde\mu)}$ to  $(Y,Z,U)$. Since $h$ is Lipschitz, we get
 \beqs
 \frac{\mathbf{E}[{K^n_T}]}{n}  = 
 \mathbf{E}\left[  \int_0^T \!\!\int_\Ic   h_i^-(s,Y^n_s,Z^n_s,U^n_s(i))\lambda(di)ds \right] \rightarrow \mathbf{E}\left[  \int_0^T \!\! \int_\Ic h_i^-(s,Y_s,Z_s,U_s(i))\lambda(di)ds \right]\;,   
  \enqs
 as $n$ goes to infinity.  Since Part (i) of Proposition \ref{MonThBSDEwJ} ensures that the sequence $(K^n_T)_{n\in\N}$ is uniformly bounded in ${\bf L^1(\Omega,\Gc_T,\P)}$, we deduce that the right hand side of the previous expression equals zero. Hence  the constraint \reff{hcons} is satisfied.\\
 
 As observed in the previous step, for any quadruple $(\check Y,\check Z,\check U,\check K)$ $\in$ ${\bf \Sc^2_{\G}}\times{\bf L^2_{\G}(W)}\times{\bf L^2(\tilde \mu)}\times{\bf A^2_{\G}}$ satisfying \reff{BSDEgen}-\reff{hcons}, the sequence $(Y^n)_{n\in\N}$ is upper bounded by $\check Y$. Passing to the limit, we deduce that $(Y,Z,U,K)$ is a minimal solution to  \reff{BSDEgen}-\reff{hcons}.\\[3mm]

\noindent {\bf Step 3: } {\it Uniqueness of the minimal solution.}\\
From the minimality condition, the uniqueness for the component $Y$ of the solution is obvious. Suppose now that we have two solutions $(Y,Z,U,K)$ and $(Y,Z',U',K')$ in ${\bf \Sc^2_{\G}}\times{\bf L^2_{\G}(W)}\times{\bf L^2(\tilde \mu)}\times{\bf A^2_{\G}}$ with $K$ and $K'$ predictable. Then we have 
\beq\label{unic1}
\int_{0}^t[f(s,Y_{s},Z_{s},U_{s})-f(s,Y_{s},Z'_{s},U'_{s})]ds+\int_{0}^t[Z'_{s}-Z_{s}]dW_{s} & &\nonumber\\ +\int_{0}^t\int_{\Ic}[U'_{s}(i)-U_{s}(i)]\mu(di,ds)+K'_{t}-K_{t} & = & 0 \;, \quad0\leq t\leq T.
\enq
Since $\mu$ is a Poisson measure, it has unaccessible jumps. Recalling that $K$ and $K'$  are predictable and taking the predictable projection in expression \reff{unic1}, we get
\beq\label{unic2}
\int_{0}^t[f(s,Y_{s},Z_{s},U_{s})-f(s,Y_{s},Z'_{s},U'_{s})]ds +\int_{0}^t[Z'_{s}-Z_{s}]dW_{s}+ K'_{t}-K_{t} & = & 0 \;, 
\enq
 for $0\le t \le T$, and 
\beqs
\int_{0}^T\int_{\Ic}[U'_{s}(i)-U_{s}(i)]\mu(di,ds) & =& 0
,\enqs
which gives $U' = U$. 
Identifying the finite variation and the Brownian parts  in \reff{unic2} we get
\beqs
\int_{0}^T[Z'_{s}-Z_{s}]dW_{s} & =& 0
,\enqs
which leads to $Z=Z'$. The uniqueness of $K$ finally follows from \reff{unic1}.
\ep

\begin{Remark}\label{RemUpperBoundY}
{\rm
Observe that the purpose of Assumption \textbf{(H1)} is simply to ensure an upper bound in ${\bf \Sc^2_{\G}}$ on the sequence of solutions $(Y^n)_{n\in\N}$ to the penalized BSDEs. 
 If such an upper bound already exists, there exists a minimal solution to  \reff{BSDEgen}-\reff{hcons} and  \textbf{(H1)} is automatically satisfied.  Hence, Theorem \ref{thmmain1} also holds under {\bf (H0)}-{\bf (H1')}.
 Particular cases where Assumption \textbf{(H1)} is satisfied are for instance presented in Theorem \ref{lingBSDECJBSDEob} below. In a Markovian setting, sufficient conditions for this assumption are also provided in Remark 3.2 of \cite{ek09}. }
\end{Remark}


\section{Connection with multidimensional reflected BSDEs}\label{SectionLinkReflBSDEs}
\setcounter{equation}{0} \setcounter{Assumption}{0}
\setcounter{Theorem}{0} \setcounter{Proposition}{0}
\setcounter{Corollary}{0} \setcounter{Lemma}{0}
\setcounter{Definition}{0} \setcounter{Remark}{0}

 In this section, we prove that one-dimensional constrained BSDEs with jumps offer a nice alternative for the representation of solutions to multidimensional reflected BSDEs studied in \cite{ht07} and \cite{hz08}. This representation has practical implications, since, for example, it opens the door to the numerical resolution of multi-dimensional reflected BSDEs via the approximation of a single one-dimensional constrained BSDE with additional artificial jumps. The arguments presented here are purely probabilistic and therefore apply in the non Markovian framework considered in \cite{ht07}. Furthermore, the proofs require precise comparison results for reflected BSDEs based on viability properties that are reported in Section \ref{App Viability}  for the convenience of the reader.

\subsection{Multidimensional reflected BSDEs}

 Recall that solving a general multidimensional reflected BSDE consists in finding $m$ triplets $(Y^{i},Z^{i},K^{i})_{i\in \Ic}$ $\in$ $({\bf \Sc_{\F}^{2}}\times{\bf L^2_{\F}(W)}\times{\bf A^2_{\F}})^m$ satisfying, for all $i\in\Ic$,
\begin{equation}\label{BSDEoblique}
\left\{
\begin{array}{l} \vspace{1mm}
Y^i_t  =  \xi^i+\int_t^T\psi_i(s,Y^1_s,\ldots,Y^m_s,Z_s^i)ds-\int_t^T\langle Z^i_s,dW_s\rangle+K^i_T-K^i_t\,, \quad 0\le t \le T\,, \\ \vspace{1mm}
Y^i_t  \geq  \max_{j\in A_{i}} h_{i,j}(t,Y^j_t)  \,, \quad 0\le t \le T\,,  \\ 
\int_0^T[Y_t^i-\max_{j\in A_{i}}\{h_{i,j}(t,Y^j_t)\}]dK^i_t=0 \;,
\end{array}
\right.
\end{equation}
where $\psi_i:~\Omega\times[0,T]\times\R^m\times\R^d\rightarrow\R$ is an $\F$-progressively measurable map, $\xi^i\in\mathbf{L}^2(\Omega,\Fc_{T}, \mathbf{P})$, $A_{i}$ is a nonempty subset of $\Ic$  
and, for any $j\in A_i \cup \{i\}$, $h_{i,j}: \Omega\times[0,T]\times\R\rightarrow \R$ is a given $\Pc_{\F}\otimes\Bc(\R)$-measurable function satisfying $h_{i,i}(t,y)=y$ for all $(t,y)\in[0,T]\times\R$. 
As detailed in Theorem 3.1 and Theorem 4.2 of \cite{hz08}, existence and uniqueness of a solution to \reff{BSDEoblique} is ensured by the following assumption:\\

\noindent \textbf{(H2)} \begin{enumerate}[(i)]
\item For any $i\in\Ic$ and $j\in A_i$, we have $\xi^{i}\geq h_{i,j}(T,\xi^{j})$.
\item For any $i\in\Ic$,  $\Eb|\xi^{i}|^2 + \Eb\int_{0}^T\sup_{y\in\R^m}|\psi_{i}(t,y,0)|^2{\bf 1}_{\{y_{i}=0\}}dt$ $<$ $+\infty$, and $\psi_i$ is  Lipschitz continuous: there exists a constant $k_{\psi}\geq 0$ such that 
\beqs
|\psi_i(t,y,z)-\psi_i(t,y',z')|\leq k_{\psi}(|y-y'|+|z-z'|)\;,\;\; (t,y,z,y',z')\in[0,T]\times[\R^m\times\R^{d}]^2\,.
\enqs
\item For any $i\in\Ic$, and $j$ $\neq$ $i$,  $\psi_i$ is nondecreasing in its $(j+1)-$th variable i.e. for any $(t,y,y',z)\in\Ic\times[\R^m]^2\times\R^d$ such that $y_{k}=y'_{k}$ for $k$ $\neq$ $j$ and $y_{j}\leq y'_{j}$, we have
\beqs 
\psi_{i}(t,y,z) & \leq & \psi_{i}(t,y',z)\quad \P-a.s.
\enqs
\item For any $(i,t,y)$ $\in$ $\Ic\times[0,T]\times\R$ and $j\in A_i$, $h_{i,j}$ is continuous,  $h_{i,j}(t,.)$ is a $1$-Lipschitz increasing function satisfying $h_{i,j}(t,y)$ $\leq$ $y$, $\P$-a.s. and we have $h_{i,j}(.,0)\in {\bf L^2(0,T)}$. 
\item For any $i\in\Ic$, $j\in A_i$ and $l\in A_j$, we have $l\in A_i\cup \{i\}$ and \beqs  h_{i,j}(t,h_{j,l}(t,y)) &<& h_{i,l}(t,y) \;, \qquad   (t,y) \in [0,T]\times\R   \;. \enqs 
\end{enumerate}

\begin{Remark}
{\rm
Part (ii) and (iii) of Assumption {\bf (H2)}  are classical Lipschitz  and monotonicity  properties of the driver. Part (iv) ensures a tractable form for the domain of $\R^m$ where $(Y^i)_{i\in\Ic}$ lies, and (i) implies that the terminal condition is indeed in the domain. Recent results in \cite{cek10} allow to relax the monotonicity condition (iii) for the case of constraint function $h$ associated to switching problems.  
}
\end{Remark}

\subsection{Corresponding constrained BSDE with jumps}

We consider now the following one-dimensional constrained BSDE with jump :  find a minimal quadruple $(\tilde Y,\tilde Z,\tilde U,\tilde K)$ $\in$ ${\bf \Sc_{\G}^2}\times{\bf L^2_{\G}(0,T)}\times{\bf
L^2(\tilde \mu)}\times{\bf A^2_{\G}}$ satisfying
\beq
\tilde Y_t &=& \xi^{I_T} + \int_t^T \psi_{I_{s^-}}(s,\tilde Y_s+\tilde U_{s}(1)\mathbf{1}_{I_{s^-}\neq 1},\ldots,\tilde Y_s+\tilde U_{s}(m)\mathbf{1}_{I_{s^-}\neq m}, \tilde Z_s) ds  + \tilde K_T - \tilde K_t  \qquad \nonumber  \\
& &\quad\quad\qquad\qquad -  \int_t^T\langle \tilde Z_s, dW_s\rangle  -  \int_t^T\int_\Ic   \tilde U_s(i)  \mu(ds,di), \;\;\;   0 \leq t \leq T,  \; a.s. ~ \label{BSDEpart}\enq
together with the constraint
\beq \label{hpart}
\mathbf{1}_{ A_{I_{t^-}}}(i)\Big[\tilde Y_{t^-}-h_{I_{t^-},i}(t,\tilde Y_{t^-}+\tilde U_t(i))\Big] & \geq & 0, \;\;\;\;\; d\Pb\otimes dt\otimes\lambda(di)
\;\;  a.e.\;,
\enq
where the process $I$ is a pure jump process defined by 
\beqs \label{dynI}
I_{t} & = & I_0 + \int_0^t \int_{\Ic}(i-I_{s^-})\mu(ds,di)\,.
\enqs

 \begin{Remark}{\rm
 If the Poisson measure rewrites $\sum_{n\geq 0}\delta_{(\kappa_{n},L_{n})}$, where $(\kappa_n)_n$ are the jump times and $(L_n)_n$ the jump sizes,  the pure jump process $I$ simply coincides with $L_{n}$ on each $[\kappa_{n},\kappa_{n+1})$.
 }\end{Remark}

 Considering $I$ as an extra source of randomness, the BSDE \reff{BSDEpart}-\reff{dynI} enters into the class of constrained BSDEs with jumps of the form \reff{BSDEgen}-\reff{hcons} studied above, with the following correspondence
 \beqs
 \xi &= &\xi^{I_{T}}\,;\\ 
 f(t,y,z,u) &=& \psi_{I_{t^-}}(t,(y+u_{i}\mathbf{1}_{I_{t^-}\neq i})_{i\in\Ic},z) \,, \qquad  (t,y,z,u)\in[0,T]\times\R\times\R^d\times\R^m     \;; \\ 
 h_i(t,y,z,v) &=& \{y-h_{I_{t^-},i}(t,y+v)\} \mathbf{1}_{i\in A_{I_{t-}}} \, \qquad(i,t,y,z,v) \in\Ic\times[0,T]\times\R\times\R^d\times\R \;.\\
 \enqs 
 
 As detailed below, Assumption \textbf{(H2)} is sufficient to ensure the existence of a one-dimensional minimal solution to the BSDE \reff{BSDEpart}-\reff{hpart}.  Remarkably, we prove hereafter that this  one-dimensional solution directly relates with the multidimensional solution of the reflected BSDE \reff{BSDEoblique}. Since the new constrained BSDE is one-dimensional, this alternative BSDE representation is promising for the numerical  resolution  of optimal switching problems . An entirely probabilistic numerical  scheme for these equations is given in \cite{ek09}.   \\
 
 We are now ready to state the main result of the paper.
 
\begin{Theorem}\label{lingBSDECJBSDEob}
 Suppose that Assumption {\bf (H2)} is in force and denote by $(Y^i,Z^i,K^i)_{i\in\Ic}$ the unique solution of \reff{BSDEoblique}. 
 Then, the constrained BSDE \reff{BSDEpart}-\reff{hpart} satisfies {\bf (H0)}-{\bf (H1)} and its unique corresponding minimal solution $(\tilde Y,\tilde Z,\tilde U,\tilde K)$ $\in$ ${\bf \Sc_{\G}^2}\times{\bf L^2_{\G}(W)}\times{\bf L^2(\tilde \mu)}\times{\bf A^2_{\G}}$ verifies
  \beq\label{linkobliqueconstrained}
 \tilde Y_t    =    Y^{I_t}_t\,, \quad   \tilde Z_t   =  Z^{I_{t^-}}_t\,, \quad \tilde U_{t}  = (Y^i_t-Y^{I_{t-}}_{t^-})_{i\in\Ic} 
 \;, \qquad 0\le t \le T \;.
 \enq
\end{Theorem}

\textbf{Proof.} The proof divides in 3 steps. First we prove the existence of a unique minimal solution to  \reff{BSDEpart}-\reff{hpart}. 
Then, we introduce a sequence of penalized BSDEs converging to the solution of the multidimensional reflected BSDE  \reff{BSDEoblique}. 
Finally, we prove that a corresponding sequence of penalized BSDEs with jumps, built via a relation of the form of \reff{linkobliqueconstrained},  converges indeed to the solution of  \reff{BSDEpart}-\reff{hpart}.\\[0mm]

\noindent {\bf Step 1: } {\it Existence and uniqueness of a minimal solution to \reff{BSDEpart}-\reff{hpart}.}\\
In order to use Theorem \ref{thmmain1}, we need to verify that Assumptions \textbf{(H0)} and \textbf{(H1)} are satisfied in this context.\\

First, parts (i) and (ii) of Assumption \textbf{(H0)} are direct consequences of  \textbf{(H2)}(ii) and (iv). 
 Fix any $(t,y,z,u,u')\in[0,T]\times\R\times\R^d\times\R^m\times\R^m$, and define $v^{(k)}\in\R^m$ by 
 \beqs
 v^{(k)} & = & (u'_{1},\ldots,u'_{k-1},u_{k},\ldots,u_{m}),\qquad 1\leq k\leq m+1\;.
 \enqs
  From the monotonicity assumption \textbf{(H2)}(iii) on the Lipschitz function $\psi$ we get
 \beqs
 f(t,y,z,u) \!-\! f(t,y,z,u') 
 \!\!\!\!& = \!\!\!\!& 
 \!\sum_{k=1}^{m}\psi_{I_{t^-}} (t, (y+v^{(k)}_i{\bf 1}_{I_{t^-}\neq i})_{i\in\Ic},z) \!-\! \psi_{I_{t^-}} (t, (y+v^{(k+1)}_i{\bf 1}_{I_{t^-}\neq i})_{i\in\Ic},z)  \\
 \!\!\!\!& \le \!\!\!\!&  k_\psi\sum_{k=1}^{m-1} (u_{k}-u'_{k})\mathbf{1}_{u_{k}\geq u'_{k}}\mathbf{1}_{k\neq I_{t^-}} \;.
 \enqs
 Taking $\gamma^{y,z,u,u'}_{t}(i)=\frac{k_{\psi}}{\lambda(i)}\mathbf{1}_{u_{k}\geq u'_{k}}\mathbf{1}_{i\neq I_{t^-}}$
 (which is well defined, since $\lambda(i)>0$ for any $i\in\Ic$), we get \textbf{(H0)}-(iii).\\
 
 In order to prove that \textbf{(H1)} holds, one needs to verify the existence of a solution to  \reff{BSDEpart}-\reff{hpart}. We indeed check hereafter that the candidate  $(\tilde Y,\tilde Z,\tilde U)$ defined in \reff{linkobliqueconstrained}  satisfies \reff{BSDEpart} as well as \reff{hpart}. Let define $N_{t}:=\mu(\Ic\times [0,t])$ for $t\in[0,T]$, the (random) number of stopping times $\kappa_{n}$, associated to the random measure $\mu$, which satisfy $\kappa_{n}\in[0,t]$.  Then,  since $Y$ is a solution of the reflected BSDE \reff{BSDEoblique}, we have
 \beqs
 Y^{L_{N_{T}}}_{\kappa_{N_{T}}} & = & \xi^{L_{N_{T}}}+\int_{\kappa_{N_{T}}}^T\psi_{L_{N_{T}}}(s,(Y^{L_{N_{T}}}_s+U_{s}(i)\mathbf{1}_{i\neq L_{N_{T}}})_{i\in\Ic},Z^{L_{N_{T}}}_s)ds\\
  & & -\int_{\kappa_{N_{T}}}^T Z^{L_{N_{T}}}_sdW_{s}+K_{T}^{L_{N_{T}}}-K_{\kappa_{N_{T}}}^{L_{N_{T}}} \;.
 \enqs
 Then, still using the equation \reff{BSDEoblique}  and identifying the jumps at  time $\kappa_{N_{T}}$, we compute: 
 \beqs
 Y^{L_{N_{T}-1}}_{\kappa_{N_{T}-1}} & = & Y^{L_{N_{T}}}_{\kappa_{N_{T}}} +\int_{\kappa_{N_{T}-1}}^{\kappa_{N_{T}}}\psi_{L_{N_{T}-1}}(s,(Y^{L_{N_{T}-1}}_s+U_{s}(i)\mathbf{1}_{i\neq L_{N_{T}-1}})_{i\in\Ic},Z^{L_{N_{T}-1}}_s)ds\\
 & & -\int_{\kappa_{N_{t}-1}}^{\kappa_{N_{t}}} Z^{L_{N_{T}-1}}_sdW_{s}+K_{\kappa_{N_{T}}}^{L_{N_{T}-1}}-K_{\kappa_{N_{T}-1}}^{L_{N_{T}-1}} +(Y^{L_{N_{T}-1}}_{\kappa_{N_{T}}} -Y^{L_{N_{T}}}_{\kappa_{N_{T}}} )\\
 & = & \xi^{I_{T}}+ \int_{\kappa_{N_{T}-1}}^T\psi_{I_{s^-}}(s,(Y^{I_{s}}_{s}+U_{s}(i)\mathbf{1}_{i\neq I_{s^-}})_{i \in \Ic},Z^{I_{s^-}}_{s})ds-\int_{\kappa_{N_{T}-1}}^T Z^{I_{s^-}}_{s}dW_{s}\\
 & & -\int_{\kappa_{N_{T}-1}}^T\int_{\Ic}U_{s}(i)\mu(di,ds)+K_{T}^{L_{N_{T}}}-K_{\kappa_{N_{T}}}^{L_{N_{T}}} +K_{\kappa_{N_{T}}}^{L_{N_{T}-1}}-K_{\kappa_{N_{T}-1}}^{L_{N_{T}-1}} \;.
 \enqs

Repeating this procedure until time $\kappa_{N_{t}+1}$ for $t\in[0,T]$, we get 
\beqs
 Y^{L_{N_{t}+1}}_{\kappa_{N_{t}+1}}  & = & \xi^{I_{T}}+ \int_{\kappa_{N_{t}+1}}^T\psi_{I_{s^-}}(s,(Y^{I_{s}}_{s}+U_{s}(i)\mathbf{1}_{i\neq I_{s^-}})_{i\in\Ic},Z^{I_{s^-}}_{s})ds-\int_{\kappa_{N_{t}+1}}^T Z^{I_{s^-}}_{s}dW_{s}\\
 & & -\int_{\kappa_{N_{t}+1}}^T\int_{\Ic}U_{s}(i)\mu(di,ds)+K_{T}^{L_{N_{T}}}-K_{\kappa_{N_{T}}}^{L_{N_{T}}} +K_{\kappa_{N_{T}}}^{L_{N_{T}-1}}-K_{\kappa_{N_{T}-1}}^{L_{N_{T}-1}}\\
& & + \ldots + K_{\kappa_{N_{t}+2}}^{L_{N_{t}+1}}-K_{\kappa_{N_{t}+1}}^{L_{N_{t}+1}}\;.
 \enqs
Combining this last expression with the equation satisfied by $Y^{L_{N_{t}}}$ between $t$ and $\kappa_{N_{t}+1}$, we deduce
the existence of a square integrable increasing process $\tilde K$ such that $(\tilde Y, \tilde Z, \tilde U, \tilde K)$ satisfies equation \reff{BSDEpart}. The reflection constraint in \reff{BSDEoblique} together with the identification  \reff{linkobliqueconstrained} imply directly that  $(\tilde Y, \tilde Z, \tilde U, \tilde K)$ satisfies the constraint \reff{hpart}. \\

 Therefore \textbf{(H0)} and \textbf{(H1)} hold for \reff{BSDEpart}-\reff{hpart} and the existence of a unique minimal solution follows from Theorem  \ref{thmmain1}.\\

 \noindent {\bf Step 2: } {\it Penalization of the multidimensional BSDE \reff{BSDEoblique}.}\\
We now introduce the following sequence of multidimensional penalized BSDEs: 
for $n\in\N$, find $m$ couples $(Y^{i,n},Z^{i,n})_{i\in \Ic}$ $\in$ $({\bf \Sc_{\F}^{2}}\times{\bf L^2_{\F}(W)})^m$ satisfying
\beq
Y^{i,n}_t  \;=\;  \xi^i+\int_t^T\psi_i^n(s,Y^{1,n}_s,\ldots,Y^{m,n}_s,Z_s^{n})ds-\int_t^T\langle Z^{i,n}_s,dW_s\rangle\;,\;\;
0 \leq t \leq T, \;\;\; i\in\Ic\;,   \quad \label{BSDEpenObl}
\enq
where the random map $\psi^n$ is defined on $[0,T]\times\R^m\times[\R^d]^m$ by
\beqs
\psi^n_i(t,y,z) &  = &  \psi_{i}(t,y,z_{i})+n\sum_{j\in A_{i}} [y_{i}-h_{i,j}(t,y_{j})]^-\lambda(j)\;, \quad (i,t,y)\in\Ic\times[0,T]\times\R^d\;.
\enqs 
For any $n\in\N$, the existence of a unique solution to \reff{BSDEpenObl} is given in the seminal paper  \cite{parpen90} 
and we prove now that the sequence of solutions to these BSDEs converges to
the solution of the multidimensional reflected BSDE \reff{BSDEoblique}.\\

 In order to prove that the sequence $(Y^{i,n})_{n\in\N}$ is nondecreasing and convergent for any $i\in\Ic$,
 we shall appeal to the multidimensional comparison theorem for reflected BSDEs presented in Section \ref{App Viability} of the paper. 
 First, since $\psi_i^n \le \psi_i^{n+1}$ for any $i\in\Ic$ and $n\in\N$, Theorem 2.1 in \cite{hp06} implies that the sequence $(Y^{.,n})_{n\in\N}$ is nondecreasing componentwise. 
 Second, we compute from the Lipschitz property of $\psi$ that
 \beqs
 - 2 \langle y , \psi^{n}(t,y',z)-\psi^{n}(t,y',z')\rangle 
 =   - 2 \langle y , \psi(t,y',z)-\psi(t,y',z')\rangle 
  \leq  k_\psi^2 \; |y|^2 + \sum_{i=1}^m |z_{i}-z_{i}'|^2 \;,
 \enqs
  $\mathbf{P}$-a.s., for any $\{t,y,y',(z,z')\}\in[0,T]\times[\R^+]^m\times\R^{m}\times[\R^{d\times m}]^2$ and $n\in\N$. 
  Therefore, since $\psi^n(t,Y_t,Z_t) = \psi(t,Y_t,Z_t)$ for $t\in[0,T]$, we deduce from Proposition \ref{thcomAp} below that 
\beq\label{majYnobl}
Y_{t}^{i,n} \leq Y_{t}^{i}\;, & \mbox{ for all } & (i,t,n)\in\Ic\times[0,T]\times \N \,.
\enq 
Introducing the sequence of processes $K^{i,n} :=  n \int_0^. \int_{A_{i}}   [Y^{i,n}_s - h_{i,j}(s,Y^{j,n}_s) ]^- \lambda(dj)ds$, for $i\in\Ic$ and $n\in\N$, we deduce from Peng's monotonic limit theorem \cite{pen99} the existence of:
\begin{itemize}
\item $\hat Y^1,\ldots,\hat Y^m$ $\F$-adapted c\`adl\`ag processes 
with $\|\hat Y^i\|_{\Sc^2} <\infty$ for all $i\in\Ic$,
\item  $\hat Z^1,\ldots,\hat Z^m$ $\in$ $\mathbf{L_{\F}^2(W)}$, 
\item $\hat K^1,\ldots,\hat K^m$ $\F$-adapted nondecreasing c\`adl\`ag processes with $\hat K_{0}^i = 0$ and $\|\hat K^i\|_{\Sc^2} <\infty$, for all $i\in\Ic$,
\end{itemize}
such that $Y^{i,n}\uparrow\hat Y^{i}$ a.e., $Y^{i,n}\rightarrow \hat Y^{i}$ in $\mathbf{L^2_{\F}(0,T)}$,  $Z^{i,n}\rightarrow \hat Z^{i}$ in $\mathbf{L^2_{\F}(W)}$ weakly, $ K ^{i,n}_{T}\rightarrow \hat K ^i_{T}$ in $\mathbf{L^2(\Omega,\Fc_{T},\P)}$ weakly and 
\begin{equation}\label{eqlim}
\left\{\begin{array}{l}
\hat Y^{i}_{t}  =  \xi^{i}+\int_t^T\psi_i(s,\hat Y^{1}_s,\ldots,\hat Y^{m}_s,\hat Z_s^{i})ds-\int_t^T\langle\hat Z^{i}_s,dW_s\rangle+\hat K^{i}_{T}-\hat K^{i}_{t}\,, \quad i\in\Ic\;, \\
\hat Y^{i}_{t}\geq\max_{j\in A_{i}}h_{i,j}(t,\hat Y^{j}_{t})\;,\quad 0\le t\le T\;,\quad i\in\Ic\;.
\end{array}\right.
\end{equation}
Observe that the last inequality  in \reff{eqlim} is not a direct consequence of Peng's monotonic limit theorem but follows instead from a similar argument as the one used in Step 2 of the proof of Theorem \ref{thmmain1} above: for $i\in\Ic$, since the sequence$(K^{i,n})_n$ is uniformly bounded in $\mathbf{L^1(\Omega,\Fc_{T},\P)}$ we have
 \beqs
 0 \; = \; \lim_{n\rightarrow\infty}  \frac{\mathbf{E}[|{K^{i,n}_T}|]}{n}  
 & = &  \lim_{n\rightarrow\infty}   \mathbf{E}\left[      \int_0^T \!\!\int_{A_{i}}   [Y^{i,n}_s - h_{i,j}(s,Y^{j,n}_s)]^- \lambda(dj)ds  \right]  \\
 &=&    \mathbf{E}\left[    \int_0^T \!\!\int_{A_{i}}   [\hat Y^{i}_s - h_{i,j}(s,\hat Y^{j}_s)]^- \lambda(dj)ds \right]\;,  \quad\qquad  i\in \Ic \;,
  \enqs
 which easily rewrites as the constraint inequality in  \reff{eqlim}. It still remains to prove that $(\hat Y, \hat Z, \hat K)$ also satisfies the  minimality property of \reff{BSDEoblique}. \\

For this purpose, we consider  the following RBSDE whose unique solution   $(\tilde Y, \tilde Z,\tilde K)$ in $({\bf \Sc_{\F}^{2}}\times{\bf L^2_{\F}(W)}\times{\bf A^2_{\F}})^m$   exists according to Theorem 2.1 in \cite{penxu04}:
\begin{equation}\label{eqqmart}
\left\{\begin{array}{l}
\tilde Y^{i}_{t}  =  \xi^{i}+\int_t^T\psi_i(s,\hat Y^{1}_s,\ldots,\hat Y^{i-1}_s,\tilde Y^{i}_s,\hat Y^{i+1}_s\ldots,\hat Y^{m}_s,\tilde Z_s^{i})ds\\ \quad \qquad\qquad \qquad \qquad \qquad \qquad \qquad-\int_t^T\langle\tilde Z^{i}_s,dW_s\rangle+\tilde K^{i}_{T}-\tilde K^{i}_{t}\,, \\ 
\tilde Y^{i}_{t}\geq\max_{j\in A_{i}}h_{i,j}(t,\hat Y^{j}_{t})\;,\quad 0\le t\le T\,, \quad i\in\Ic\;, \\
\int_{0}^T [\tilde Y ^i_{t^-}-\max_{j\in A_{i}}h_{i,j}(t,\hat Y_{t^-}^j)]d\tilde K_{t}^j=0\;,\quad i\in\Ic\;.
\end{array}\right.
\end{equation}

We note that \reff{eqlim} and \reff{eqqmart} have the same lower barrier. For any $i\in\Ic$, since $\tilde Y ^i$ is the smallest $\psi_{i}$-supermartingale with lower barrier $\max_{j\in A_{i}} h(., \hat Y ^j_{.})$, we know from Theorem 2.1 in  \cite{penxu04} that  $\tilde Y  ^i\leq \hat Y  ^i$. 
  
On the other hand, we deduce from \textbf{(H2)} (iii) that
\beqs
\psi_{i}^n(s,\hat Y ^1_{s},\ldots,\hat Y ^{i-1}_{s},y,\hat Y ^{i+1}_{s}, \ldots, \hat Y ^{m}_{s}) & \geq & \psi_{i}^n(s, Y ^{1,n}_{s},\ldots, Y ^{i-1,n}_{s},y, Y ^{i+1,n}_{s}, \ldots, Y ^{m,n}_{s})\;,
\enqs
for all $(i,s,y,n)\in\Ic\times[0,T]\times\R\times\N$, $\P$-a.s..
 For $i\in\Ic$, since $\tilde Y ^i \geq \max_{j\in A_{i}}h_{i,j}(.,Y^j_{.})$, combining \textbf{(H2)} (iv) and a comparison theorem for one dimensional reflected BSDEs, we get $Y^{i,n}\leq \tilde Y^i$ for any $n\in\N$, and, sending $n$ to infinity, deduce $\hat Y^{i}\leq \tilde Y^i$. 
 
 Therefore $\hat Y = \tilde Y$ and  $(\hat Y,\hat Z,\hat K)$ satisfies
\begin{equation}\label{eqreffinale}
\left\{\begin{array}{l}
\hat Y^{i}_{t}  =  \xi^{i}+\int_t^T\psi_i(s,\hat Y_s,\hat Z_s^{i})ds
-\int_t^T\langle\hat Z^{i}_s,dW_s\rangle+\hat K^{i}_{T}-\hat K^{i}_{t} \;,\quad 0\le t\le T\,, \quad i\in\Ic\;, \\
\hat Y^{i}_{t}\geq\max_{j\in A_{i}}h_{i,j}(t,\hat Y^{j}_{t})\;,\quad 0\le t\le T\,, \quad i\in\Ic\;, \\
\int_{0}^T [\hat Y ^i_{t^-}-\max_{j\in A_{i}}h_{i,j}(t,\hat Y_{t^-}^j)]d\hat K_{t}^j=0\;,\quad i\in\Ic\;.
\end{array}\right.
\end{equation}
Notice that the minimality condition in \reff{eqreffinale} differs from the expected one in \reff{BSDEoblique}. Nevertheless, those two coincide whenever $\hat Y$ is continuous, property that we verify now.\\

Suppose on the contrary that    $ \hat Y ^{i_{1}}_{t}\neq \hat Y ^{i_{1}}_{t^-}$ for some fixed $(i_{1},t)\in\Ic\times[0,T]$. Then, we deduce from \reff{eqreffinale} that $ \hat Y ^{i_{1}}_{t} - \hat Y ^{i_{1}}_{t^-}= \hat K ^{i_{1}}_{t^-} - \hat K ^{i_{1}}_{t}< 0$,  which further implies 
\beqs
\hat Y ^{i_{1}}_{t^-}=\max_{j\in A_{i}}h_{i_{1},j}(t,\hat Y_{t^-}^j)=h_{i_{1},i_{2}}(t,\hat Y_{t^-}^{i_{2}})\;,
\enqs
for some $i_{2}\neq i_{1}$. Using the constraint satisfied by $\hat Y$, we get
\beqs
h_{i_{1} ,i_{2}} (t, \hat Y^{i_{2}}_{t^- }) = \hat  Y^{i_{1}}_{t^-}>  \hat Y^{i_{1}}_{t} \geq \max_{i\in A_{i_{1}}}  
h_{i_{1} ,i} (t, \hat  Y ^i_{t} ) \geq h_{i_{1} ,i_{2}} (t, \hat  Y^{i_{2}}_{t} ). 
\enqs 
Thus $ \hat Y^{i_{2}}_{t} < \hat Y^{i_{2}}_{t^-}$. Repeating this argument we get a finite cyclic sequence $(i_{k})_{1\leq k\leq N}$ such that $i_{N}=i_{1}$ and 
\beqs
\hat Y^{i_{k-1}} _{t^-}  & = &  h_{i_{k-1},i_{k}}(t, \hat Y^{i_{k}}_{t^-} )\;,\quad 2\leq k \leq N\;, 
\enqs
which contradicts \textbf{(H2)} (v).\\[2mm]

 \noindent {\bf Step 3: } {\it Link between solutions of  BSDE \reff{BSDEoblique} and BSDE  \reff{BSDEpart}-\reff{hpart}.}\\
 For $n\in\N$, define the process $(Y^{I,n},Z^{I,n},U^{I,n})$ $\in$   ${\bf \Sc_{\G}^2}\times{\bf L^2_{\G}(W)}\times{\bf
L^2(\tilde \mu)}$  by
\begin{align}\label{defYi}
Y^{I,n}_t \,  := \,  Y^{I_t,n}_t\,,\;\;  Z^{I,n}_t  \,:= \, Z^{I_{t^-},n}_t 
  \;\;\mbox{ and } \;\; U^{I,n}_t\,:=\, (Y^{i,n}_t-Y^{I,n}_{t^-} )_{i\in\Ic}  
  \,,   \quad  0 \le t \le T.
\end{align}
 In order to obtain the correspondence  \reff{linkobliqueconstrained}, it only remains to prove that $(Y^{I,n},Z^{I,n},U^{I,n})_n$ converges to $(\tilde Y,\tilde Z,\tilde U)$.\\

As in Step 1, writing  the dynamics of \reff{BSDEpenObl} between each successive stopping times  associated to the random measure $\mu$, we easily check that $(Y^{I,n},Z^{I,n},U^{I,n})$ is the unique solution of the following penalized BSDE 
\beq
Y^{I,n}_t
 &\!\!\! = \!\!& \xi^{I_T}+\int_t^T\psi_{I_{s^-}}(s,Y^{I,n}_{s}+U^{I,n}_s(1)\mathbf{1}_{I_{s^-}\neq 1},\ldots,Y^{I,n}_{s}+U^{I,n}_s(m)\mathbf{1}_{I_{s^-}\neq m},Z_s^{I,n})ds \nonumber \\
  &\!\!\! - \!\!& \int_t^T\langle Z^{I,n}_s,dW_s\rangle+n\int_t^T\hspace{-2mm}\int_{\Ic}h_{i}^-(s,Y^{I,n}_{s^-},Z^{I,n}_{s},U^{I,n}_s(i))\lambda(di)ds+\int_{t}^T\hspace{-2mm}\int_{\Ic}U^{I,n}_{s}(i)\mu(ds,di)
, \nonumber
\enq
 for $0 \le t \le T$.
 Therefore, Step 1 
  ensures that we can apply Theorem \ref{thmmain1} and
 we get
 \beq\label{cvlink} 
 \|Y^{I,n}-\tilde  Y\|_{_{{\bf L^2(0,T)}}} + \|Z^{I,n} -  \tilde Z\|_{_{{\bf L^p(0,T)}}} + \|U^{I,n} - \tilde U\|_{_{{\bf L^p(\tilde\mu)}}}   & \longrightarrow & 0\,,\quad p<2,
 \enq 
 where we recall that $(\tilde Y,\tilde Z,\tilde U)$ is the minimal solution to \reff{BSDEpart}-\reff{hpart}.
 Combining this result with \reff{defYi} and Step 2  concludes the proof.
\ep
\\


\section{Subsidiary technical points}\label{App Comp Jumps}\label{SecAppendix}
\setcounter{equation}{0} \setcounter{Assumption}{0}
\setcounter{Theorem}{0} \setcounter{Proposition}{0}
\setcounter{Corollary}{0} \setcounter{Lemma}{0}
\setcounter{Definition}{0} \setcounter{Remark}{0}

This section regroups technical properties which are mainly extensions of existing results but that we could not find as such in the literature. They are not the main focus of the paper but still present some interest in themselves. This dissociation allows to present them in a more abstract setting and simplifies their possible future quotation.  We provide a comparison and a monotonic limit theorem for BSDEs with jumps, as well as viability and comparison properties for multidimensional reflected BSDEs.

\subsection{A comparison theorem for reflected BSDEs with jumps}\label{App Comp Jumps}

 We derive here a general comparison theorem for reflected BSDEs with jumps. This extends the results  of Theorem 2.5 in \cite{roy06} obtained in the non-reflected case.

\begin{Proposition}\label{PropComparisonBSDE}
Let $f_{1},f_{2}: \Omega\times[0,T]\times\R\times\R^d\times\R^m\rightarrow \R$ two generators satisfying Assumption \textbf{(H0)} and $\xi_{1},\xi_{2}$ $\in$ $\mathbf{L}^2(\Omega,\Gc_{T},\mathbf{P})$.   Let $(Y^1,Z^{1},U^1)\in\Sc^2_{\G}\times\mathbf{L^2_{\G}(W)}\times\mathbf{L}^2(\tilde \mu)$ satisfying on $[0,T]$
\beq
Y^1_{t} & = & \xi^1+\int_{t}^Tf_{1}(s,Y^1_{s},Z^1_{s},U^1_{s})ds-\int_{t}^T\langle Z^1_{s},dW_{s}\rangle-\int_{t}^T\int_{\Ic}U^1_{s}(i)\mu(ds,di) \;,
\enq
and $(Y^2,Z^{2},U^2,K^2)\in\Sc^2_{\G}\times\mathbf{L^2_{\G}(W)}\times\mathbf{L}^2(\tilde \mu)\times\mathbf{A^2_{\G}}$ satisfying on $[0,T]$
\begin{align}
Y^2_{t} & = & \xi^2+\int_{t}^Tf_{2}(s,Y^2_{s},Z^2_{s},U^2_{s})ds-\int_{t}^T\langle Z^2_{s},dW_{s}\rangle-\int_{t}^T\int_{\Ic}U^2_{s}(i)\mu(ds,di)
+K^2_{T}-K^2_{t}\;.\qquad~
\end{align}
If $\xi^1$ $\leq$ $\xi^2$ and $f_{1}(t,Y^1_{t},Z^1_{t},U^1_{t})$ $\leq$ $f_{2}(t,Y^1_{t},Z^1_{t},U^1_{t})$ for all $t$ $\in$ $[0,T]$, then we have 
 \beqs
 Y^1_{t} & \leq & Y^2_{t}\;, \qquad 0\le t\le T \;.
 \enqs
\end{Proposition}

\ni\textbf{Proof.}
Let us denote $\bar{Y}$ $:=$ $Y^2-Y^1$,  $\bar{Z}$ $:=$ $Z^2-Z^1$, $\bar U$ $:=$ $U^2- U^1$, $\bar{f}$ $=$ $f_{2}(.,Y^2,Z^2,U^2)-f_{1}(.,Y^1,Z^1,U^1)$ 
and $\bar \xi$ $=$ $\xi^2-\xi^1$ 
 so that 
 \begin{align}\label{EDSRYbar}
 \bar{Y}_t &=&\bar \xi + \int_t^T \bar f_{s}ds-\int_t^T \langle\bar{Z}_{s},dW_{s}\rangle -\int_t^T\hspace{-2pt}\int_\Ic\bar{U}_s(i)\mu(ds,di) 
+K^2_T - K^2_t ,\;\; 0\le t \le T\;, \qquad~
\end{align}
Let now define the process $a$ by
 \beqs 
 a_t &:=& \frac {f_{2}(t,Y^2_t,Z^2_t,U^2_t)-f_{2}(t,Y^1_t,Z^2_t,U^2_t)}{\bar{Y}_t}    \mathbf{1}_{\{\bar{Y}_{t}\neq 0\}}\,, \quad 0\le t\le T\;,
 \enqs
 and $b$ the $\mathbb{R}^d$-valued process defined component by component by
 \beqs 
 b_t^k &:=&
 \frac{f_{2}(t,Y^1_t,Z^{(k-1)}_t,U^2_t)-f_{2}(t,Y^1_t,Z^{(k)}_t,U^2_t)}{V_t^{k}}
 \mathbf{1}_{\{V_t^{k}\neq 0\}}
 \,,\;\; k=1,\ldots,d\;,   \quad 0\le t\le T\;,
 \enqs
where $Z_t^{(k)}$ is the $\mathbb{R}^d$-valued random vector whose $k$ first components are those of  $Z^1$ and whose  $(d-k)$ lasts are those of  $Z^2$, and  $V_t^{k}$ is the  $k$-th component of $Z^{(k-1)}_t-Z^{(k)}_t$.
 
Notice that the processes $a$ and $b$ are $\P$-$a.s.$ bounded
since $f_2$ is Lipschitz continuous. 
Observe also that the process $\bar  K$ defined on $[0,T]$ by  
\beqs
\bar K_t \, := \,  K^2_{t}-\int_{0}^t\!\!\int_{\Ic}{^2\gamma_{s}^{Y_{s^-}^1,Z^1_{s},U^1_{s},U^2_{s}}}\bar{U}_{s}(i)\lambda(di)ds+\int_{0}^t\!(f_{2}(s,Y^1_{s},Z^{1}_s, U^2_s)-f_{1}(s,Y^1_{s},Z^{1}_s,U^{1}_s))ds
\enqs 
 is a non-decreasing process since $f_2$ satisfies \textbf{(H0)} (iii) with associated bounded process  $^2\gamma$, and $f_{1}(t,Y^1_{t},Z^1_{t},U^1_{t})$ $\leq$ $f_{2}(t,Y^1_{t},Z^1_{t},U^1_{t})$, for all $t$ $\in$ $[0,T]$.
 With these notations,  
we rewrite \reff{EDSRYbar} as:
 \beqs \bar{Y}_t &=& \bar \xi + \int_t^T\left(a_s \bar{Y}_s+\langle b_s ,\bar{Z}_s\rangle +\int_{\Ic} 
  {^2\gamma_{s}^{Y^1_{s^-},Z^1_{s},U^1_{s},U^2_{s}}}(i)
  \bar{U}_{s}(i) \lambda(di)\right)ds  \\
  & & 
\qquad\qquad\qquad\qquad -\int_t^T\langle\bar{Z}_s, dW_s\rangle
- \int_t^T\int_\Ic \bar U_s(i) \mu(ds,di) +\bar{K}_T - \bar{K}_t\,. 
\enqs Consider now the positive process  $\Gamma$  solution of
the s.d.e.:
 \beqs d\Gamma_t & = & \Gamma_{t^-} \left(a_t dt+ \langle b_{t} , dW_t\rangle+\int_{\Ic}  {^2\gamma_{t}^{Y^1_{t^-},Z^1_{t},U^1_{t},U^2_{t}}} (i)\mu(dt,di)\right), \;\;\; \Gamma_0 \; = \;  1.
\enqs 
Notice that $\Gamma$ lies in $\Sc^2_{\G}$ since $a$, $b$  and $\gamma$ are bounded, and $\Gamma$ is positive since $^2\gamma>-1$.
 A direct application of It\^o's formula leads to 
\beqs
d(\Gamma\bar{Y})_t & = & \langle\Gamma_{t^-}\bar Z_{t} + \bar Y_{t^-}\Gamma_{t^-}b_t,dW_{t}\rangle+\Gamma_{t^-} \int_{\Ic}{^2\gamma_{t}^{Y^1_{t^-},Z^1_{t},U^1_{t},U^2_{t}}}(i)\bar U_{t}(i)\tilde \mu(ds,di)-\Gamma_{t^-}d\bar K_{t}\,,
\enqs 
recall that $\tilde \mu$ is the compensated measure associated to $\mu$. 
 Hence,  the process $\Gamma\bar Y$ is a supermartingale since $\Gamma$ $>$  $0$.  Therefore
\beqs
 \Gamma_t \bar{Y}_t
 & \geq & \E \left[\left.\Gamma_T \bar{Y}_T\right|\Gc_{t}\right]~ =~ \E \left[\left.\Gamma_T \bar{\xi}\right|\Gc_{t}\right]\geq 0\,, \qquad 0\leq t\leq T\,,
\enqs
leading to $\bar{Y}\geq 0$.
\ep

\subsection{Monotonic limit theorem for BSDE with jumps}\label{App Monotonic Jumps}

This paragraph is devoted to the extension of Peng's monotonic limit theorem to the framework of BSDEs driven by a Brownian motion and a Poisson random measure. In the particular case where the driver $f$ does not depend on the jump component $U$, this extension can be obtained combining several results derived in Section 3 of \cite{kmpz08}. For sake of completeness, we provide here a proof of the result. 

\vspace{2mm}

We consider a sequence $(Y^n,Z^n,U^n,K^n)_n$ in $\mathbf{\Sc^2_\G}\times \mathbf{L^2_\G(W)}\times \mathbf{L^2(\tilde \mu)}\times\mathbf{A^2_\G}$ such that 
\beqs
Y^n_{t} \, = \, Y^n_{T}+\int_{t}^Tg(s,Y^n_{s},Z^n_{s},U^n_{s})ds-\int_{t}^T\langle Z^n_{s},dW_{s}\rangle-\int_{t}^T\int_{\Ic}U^n_{s}(i)\mu(ds,di)+K^n_{T}-K^n_{t}\;,
\enqs 
for all $t\in[0,T]$ and all $n\in\N$. Here  $g$ $:$ $\Omega\times[0,T]\times\R\times\R^d\times\R^m \rightarrow\R$, is $\mathfrak{P}_{\G}\otimes \Bc(\R)\otimes \Bc(\R^d)\otimes \Bc(\R^m)$-measurable, 
 We also introduce the following assumption :

\newpage

\ni\textbf{(H3)}\begin{enumerate}[(i)]

\item $g(.,0,0,0)$ is square integrable:
\beqs
\E\int_0^T|g(t,0,0,0)|^2dt & < & \infty\;.
\enqs

\item There exists a constant $k>0$ such that the function $g$  satisfies  $\P$-a.s. the uniform Lipschitz property:
\beqs
 |g(t,y,z,u)-g(t,y',z',u')|  & \leq & k |(y,z,u)-(y',z',u')|\;, 
 \enqs 
  for all  
 $\{t,(y,z,u),(y',z',u')\}\in[0,T]\times[\R\times\R^{d}\times\R^m]^2$.

\item For any $t\in[0,T]$, $(Y^n_t)_n$ converges increasingly to $Y_t$, and we have  $\|Y\|_{\Sc^2}< \infty$.

\item $K^n$ is a continuous process, for any $n\in\N$.

\end{enumerate}

\vspace{2mm}

\begin{Proposition}\label{MonThBSDEwJ}

Suppose that \textbf{(H3)} holds true. Then, we have:
 \begin{enumerate}[(i)]
 \item Up to a modification, $Y\in{\bf \Sc^2_{\G}}$ and there exists $(Z,U,K)\in {\bf L^2_{\G}(W)}\times{\bf L^2(\tilde\mu)}\times{\bf A^2_{\G}}$ with $K$ $\G$-predictable, such that
\beqs 
 \|Y^n-Y\|_{_{{\bf L^2(0,T)}}} + 
\|Z^n-Z\|_{_{{\bf L^p(0,T)}}} + \|U^n-U\|_{_{{\bf L^p(\tilde\mu)}}} 
   & \longrightarrow_{n \rightarrow\infty} & 0,  \qquad  1\le p < 2 \;,
 \enqs 
 and $K_{t}$ is the weak limit of $(K^n_{t})_{n\in\N}$ in ${\bf L^2(\Omega,\Gc_{t},\P)}$, for any $t\in[0,T]$. 
 Moreover, $(Z,U,K)$ is the weak limit of $(Z^n,U^n,K^n)_{n\in\N}$ in ${\bf L^2_{\G}(W)}\times{\bf L^2(\tilde\mu)}\times{\bf L^2_{\G}(0,T)}$. 
\item The quadruple $(Y,Z,U,K)$ satisfies
\beq\nonumber 
Y_{t} & = & Y_{T}+\int_{t}^Tg(s,Y_{s},Z_{s},U_{s})ds-\int_{t}^T\langle Z_{s},dW_{s}\rangle\\
 & & \qquad \qquad \qquad-\int_{t}^T\int_{\Ic}U_{s}(i)\mu(ds,di)+K_{T}-K_{t}\;,\qquad 0\leq t\leq T \,.\label{eqmonotonic}
\enq 
 \end{enumerate}
  \end{Proposition}

\ni\textbf{Proof.} The proof of Proposition \ref{MonThBSDEwJ} is an adaptation of the proof of Theorem 2.4 in \cite{pen99}. The main assumption which allows to extend  the arguments of  \cite{pen99} is the continuity of each process $K^n$, $n\in\N$. We recall the main steps of the proof and explain how the continuity assumption provides the result.

\paragraph{1.}\textit{Uniform estimate.} Since the sequence $(Y^n)_n$ is monotonic, there exists a constant $C$ such that 
\beq\label{Eq S2 Yn}
\sup_{n\in\N}\|Y^n\|_{\Sc^2} & \leq & \|Y^0\|_{\Sc^2}  +  \|Y\|_{\Sc^2} \;\le\; C\;.
\enq
 Applying It\^o's formula to $|Y^n|^2$ and using \textbf{(H3)} (ii), we have
 \beqs
 \E|Y^n_t|^2 & = & \E|Y^n_T|^2 + 2\E\int_t^T\!\! Y_s^ng(s,Y_s^n,Z_s^n,U_s^n)ds-\E\int_t^T|Z^n_s|^2ds
\\ && -\E\int_t^T\!\!\!\int_\Ic\!\big(|Y^n_{s^-}+U_s^n(i)|^2-|Y^n_{s^-}|^2\big)\mu(di,ds)+2\E\int_t^TY_s^ndK_s^n\\
 & \leq & \E|Y^n_T|^2 + 2\E\int_t^T\!\! |Y_s^n|\big(g(s,0,0,0)+k|Y_s^n|+k|Z_s^n|+k|U_s^n|\big)ds-\E\int_t^T|Z_s^n|ds
\\ &&-\E\int_t^T\!\!\!\int_\Ic\!\big(2Y^n_{s}U_s^n(i)-|U^n_{s}(i)|^2\big)\lambda(di)ds +2\E\sup_{s\in[0,T]}|Y_s^n|\int_t^TdK_s^n \;.\\
 \enqs
 Using the inequality $2ab\leq \eta|a|^2+{|b|^2\over\eta}$ for $a,b\in\R$ and $\eta>0$ and \textbf{(H3)} (i), we get the existence of a constant $C$ s.t.
 \beq
 \E\int_0^T|Z_s^n|^2ds+\E\int_0^T\int_\Ic|U_s^n(i)|^2\lambda(di)ds  \leq   C\Big(\E\sup_{t\in[0,T]}|Y^n_t|^2+1\Big) 
+ 2\E K^n_T\sup_{t\in[0,T]}|Y^n_t| \;.  \quad\label{eqUnifest1}
 \enq
  Then since 
  \beqs
  K^n_T & = & Y_0^n -Y_T^n -\int_0^Tg(s,Y_s^n,Z^n_s,U_s^n)ds +\int_0^T\langle Z^n_s, dW_s \rangle+\int_0^T\int_\Ic U^n_s(i)\mu(di,ds)\;,
  \enqs
  we have from (H3) (ii), the existence of a positive constant $C'$ s.t.
  \beq\label{eqUEK}
  \E | K^n_T |^2 \, \leq \, C'\Big(1+\E\sup_{t\in[0,T]} |Y^n_t|^2 + \E\int_0^T |Z_t^n|^2dt+\E\int_0^T\int_\Ic|U^n_s(i)|^2\lambda(di)ds\Big) \;.
  \enq
Applying the inequality $2ab\leq 2C'|a|^2+{|b|^2\over{2C'}}$ for $a,b\in\R$, we obtain
\beqs
2\E K^n_T\sup_{t\in[0,T]}|Y^n_t| \, \leq \, {1\over 2}\E\int_0^T|Z_s^n|^2ds+ {1\over 2}\E\int_0^T\int_\Ic|U_s^n(i)|^2\lambda(di)ds 
+ C''\Big(1+\E\sup_{t\in[0,T]}|Y^n_t|^2\Big)\;.
\enqs
Combining this last estimate with \reff{Eq S2 Yn} and \reff{eqUnifest1}, we obtain  a constant $C$ such that
\beqs
\|Y^n\|_{{\bf \Sc^2}}+\|Z^n\|_{{\bf L^2(0,T)}}+\|U^n\|_{{\bf L^2(\tilde\mu)}} & \leq & C\;, \qquad n\in\N\,.
\enqs
Then combining the previous inequality with \reff{eqUEK} we get 
\beq\label{UE2}
\|Y^n\|_{{\bf \Sc^2_{\G}}}+\|Z^n\|_{{\bf L^2(0,T)}}+\|U^n\|_{{\bf L^2(\tilde\mu)}}+\|K^n\|_{{\bf \Sc^2_{\G}}}& \leq &C \;, \qquad n\in\N\,.
\enq 

\paragraph{2.}\textit{Weak convergence.} Using the previous uniform estimate and the Hilbert structure of ${\bf L^2_{\G}(W)}\times{\bf L^2(\tilde\mu)}\times {\bf L^2_{\G}(0,T)}\times \bf{L^2_{\G}(0,T)}$, we  deduce the existence of a subsequence of $(Z^n,U^n,K^n, g(.,Y^n,Z^n,U^n))_{n}$, which converges weakly to some process $(Z,U,K,G)$ in ${\bf L^2_{\G}(W)}\times{\bf L^2(\tilde\mu)}\times {\bf L^2_{\G}(0,T)}\times \bf{L^2_{\G}(0,T)}$. 

Identifying the limits of $(Y^n)_{n}$ and  $(Z^n,U^n,K^n, g(.,Y^n,Z^n,U^n))_{n}$, we get 
\beq\label{identifLim}
Y_{t}  \,=\, Y_{T}+ \int_{t}^TG_{s}ds-\int_{t}^T\langle Z_{s},dW_{s}\rangle-\int_{t}^T\int_{\Ic}U_{s}(i)\mu(ds,di) +K_{T} - K_{t},\quad 0\leq t \leq T.\quad
\enq
The predictability of the process $K$ comes from the predictability of each $K^n$ and the completeness of $\mathbf{L^2_{\G}(W)}$ for the weak topology. 
\paragraph{3.}\textit{Properties of the process} $K.$ We first observe from Lemma 2.2 in \cite{pen99} that the process $K$ admits a c\`adl\`ag modification. We then establish that the contribution of the jumps of $K$ is mainly concentrated within a finite number of intervals with sufficiently small total length. 

As in Lemma 2.3 in \cite{pen99}, 
for any $\delta, \epsilon>0$, there exists a finite number of pairs of stopping times $(\sigma_{k},\tau_{k})_{0\leq k \leq N}$ with $0<\sigma_{k}\leq \tau_{k}\leq T$ such that
\begin{enumerate}[(i)]
\item $(\sigma_{j},\tau_{j}]\cap(\sigma_{k},\tau_{k}]=\emptyset$ for $j\neq k$;
\item $\Eb\sum_{k=0}^N(\tau_{k}-\sigma_{k})\geq T-{\eps}$;
\item $\Eb\sum_{k=0}^N\sum_{\sigma_{k}<t\leq \tau_{k}}|\Delta K_{t}|^2\leq \delta$.
\end{enumerate}
This result is derived with similar arguments as in \cite{pen99}, relying only on the right continuity of the filtration and the predictability of the process $K$. 
More precisely, its proof is based on Lemma A.1 in \cite{pen99} and the fact that, since $K$ is predictable, its jump times are predictable stopping times and hence could be announced. In other words, if $\tau$ is a jump time of $K$ then there exist a sequence of stopping times $(\tau_{k})_{k}$ with $\tau_{k}<\tau$ for each $k$ and $\tau_{k}\uparrow\tau$ as $k$ goes to infinity.  Combining these two results, we end this step as in \cite{pen99}.  
\paragraph{4.}\textit{Strong convergence.}
From the previous step, for any $\delta,~\eps>0$, there exists a finite number of disjoint stochastic intervals $(\sigma_k,\tau_k]$, $k=0,\ldots,N$,    satisfying
\begin{enumerate}[(i)]
\item $\Eb\sum_{k=0}^N(\tau_{k}-\sigma_{k})\geq T-{\eps\over 2}$;
\item $\Eb\sum_{k=0}^N\sum_{\sigma_{k}<t\leq \tau_{k}}|\Delta K_{t}|^2\leq {\delta\eps\over 3}$.
\end{enumerate}

Then applying It\^o's formula to $|Y^n-Y|$ on $(\sigma_k,\tau_k]$ and summing over $k$ we have
\beqs
 \E\int_{\sigma_k}^{\tau_k}|Z^n_s-Z_s|^2ds + \E\int_{\sigma_k}^{\tau_k}\int_\Ic |U_s(i)-U_s^n(i)|^2\lambda(di)ds  \leq  \qquad\qquad\qquad\qquad\qquad\qquad \\ 
\qquad |Y^n_{\tau_k}-Y_{\tau_k}|^2 + \sum_{t\in(\sigma_k,\tau_k]}|\Delta K_t|^2 \label{encore}
+ 2 \E\int_{\sigma_k}^{\tau_k}  |Y^n_{s}-Y_{s}||g(s,Y^n_s,Z^n_s,U^n_s)-G_s|ds  \\
+2 \E\int_{\sigma_k}^{\tau_k}\int_\Ic\big(|U_s(i)-U_s^n(i)||Y_s-Y_s^n|\big)\lambda(di)ds\;
+2 \E\int_{\sigma_k}^{\tau_k} (Y_s-Y^n_s) dK_s  \; .  \nonumber
\enqs
Hence, summing over $k$, we obtain the existence of a constant $C$ such that
\beq\nonumber
\sum_{k=0}^N \E\int_{\sigma_k}^{\tau_k}|Z^n_s-Z_s|^2ds+\E\int_{\sigma_k}^{\tau_k}\int_\Ic |U_s(i)-U_s^n(i)|^2\lambda(di)ds  \leq  \qquad\qquad\qquad\qquad\qquad\\ 
 C\Big( \E\int_0^T |Y_s-Y^n_s| dK_s +\E\int_0^T|Y_{s}-Y^n_{s}|(|g(s,Y^n_s,Z^n_s,U^n_s)-G_s|+1)ds\Big) \nonumber \\
+\; C\Big(\sum_{k=0}^N\E|Y^n_{\tau_k}-Y_{\tau_k}|^2+\E\sum_{k=0}^N\sum_{t\in(\sigma_k,\tau_k]}  |\Delta K_t|^2  \Big) \,.\qquad\qquad\qquad\qquad\qquad\qquad\label{encore}
\enq
Using Cauchy Schwartz inequality, we have
\beq\label{cvzero1}
\E\int_0^T|Y^n_{s}-Y_{s}||g(s,Y^n_s,Z^n_s,U^n_s)-G_s|ds & \leq & C\Big(\E\int_0^T|Y^n_{s}-Y_{s}|^2ds\Big)^{1\over 2}~\longrightarrow~ 0\;,\qquad 
\enq
as $n\rightarrow\infty$. Moreover, since $|Y^0_s-Y_s|\geq |Y^n_s-Y_s|\rightarrow 0$ and
\beqs
\E\int_0^T|Y^0_s-Y_s|dK_s & \leq & \Big(\E\sup_{[0,T]}|Y^0-Y|^2\Big)^{1\over 2}\big(\E|K_T|^2\big)^{1\over 2}  <~\infty\;,
\enqs
we get from the dominated convergence theorem that
\beq\label{cvzero2}
\E\int_0^T|Y^n_s-Y_s|dK_s & \longrightarrow & 0 \;, \qquad \mbox{ as $n\rightarrow\infty$}\;.
\enq
Finally, since 
\beqs
\sum_{k=0}^N\E|Y^n_{\tau_k}-Y_{\tau_k}|^2 & \leq & N\E\sup_{[0,T]}|Y^0-Y|^2 ~<~\infty\;,
\enqs
we get from the dominated convergence theorem that
\beq\label{cvzero3}
\sum_{k=0}^N\E|Y^n_{\tau_k}-Y_{\tau_k}|^2 & \longrightarrow & 0 \;, \qquad \mbox{ as $n\rightarrow\infty$}\;.
\enq
Combining \reff{encore} with \reff{cvzero1}, \reff{cvzero2} and \reff{cvzero3}, we get
\beqs
\overline{\lim_{n\to\infty}}\sum_{k=0}^N \E\int_{\sigma_k}^{\tau_k}|Z^n_s-Z_s|^2ds+\E\int_{\sigma_k}^{\tau_k}\!\!\!\int_\Ic |U_s(i)-U_s^n(i)|^2\lambda(di)ds  \leq \sum_{k=0}^N\sum_{t\in(\sigma_k,\tau_k]}\hspace{-2mm}|\Delta K_t|^2 \leq \frac{\eps\delta}{3}\;.
\enqs 
Thus, there exists an integer $l_{\eps,\delta}$ such that 
\beqs
\sum_{k=0}^N \E\int_{\sigma_k}^{\tau_k}|Z^n_s-Z_s|^2ds+\E\int_{\sigma_k}^{\tau_k}\int_\Ic |U_s(i)-U_s^n(i)|^2\lambda(di)ds  & \leq & \frac{\eps\delta}{2}\;,\qquad \mbox{for any $n\geq l_{\eps,\delta}\;$.}
\enqs
Therefore, in the product space $([0,T]\times \Omega,\Bc([0,T])\otimes\Gc)$,  we have
\beqs
m \otimes \P\Big((s,\omega)\in\cup_{k=0}^N(\sigma_k,\tau_k]\times\Omega, ~|Z^n_s-Z_s|^2\geq \delta\Big) & \leq & {\eps\over 2}\;,
\enqs
 and, in the product space $([0,T]\times \Omega\times \Ic,\Bc([0,T])\otimes\Gc\otimes\sigma(\Ic))$  we have 
\beqs
m \otimes \lambda \otimes \P\Big((s,i,\omega)\in\cup_{k=0}^N(\sigma_k,\tau_k]\times\Ic\times\Omega, ~|U^n_s(i)-U_s(i)|^2\geq \delta\Big) & \leq & {\eps\over 2}\;,
\enqs
where $m$ denotes the lebesgue measure on $\R^+$.  This implies that 
\beqs
\lim_{n\to\infty} m\otimes \P\Big((s,\omega)\in\cup_{k=0}^N(\sigma_k,\tau_k]\times\Omega, ~|Z^n_s-Z_s|^2\geq \delta\Big) & = & 0\;.
\enqs
and 
\beqs
\lim_{n\to\infty} m\otimes \lambda \otimes \P\Big((s,i,\omega)\in\cup_{k=0}^N(\sigma_k,\tau_k]\times\Ic\times\Omega, ~|U^n_s(i)-U_s(i)|^2\geq \delta\Big) & = & 0\;.
\enqs
Hence $(Z^n)_n$ (resp. $(U^n)_n$) converges in measure to $Z$ (resp. $U$) and since it is bounded in $\mathbf{L^2(0,T)}$ (resp.  $\mathbf{L^2(\tilde \mu)}$), it converges  in $\mathbf{L^p(0,T)}$ (resp. $\mathbf{L^p(\tilde \mu)}$) for all $p<2$.
Then, combining \textbf{(H3)} (i) with the previous strong convergence of $(Z^n,U^n)_n$ to $(Z,U)$ we get
\beqs
G_s & = & g(s,Y_s,Z_s,U_s)\;, \qquad 0\leq s \leq T\;,
\enqs   
and from \reff{identifLim}, we deduce that $(Y,Z,U,K)$ satisfies \reff{eqmonotonic}.
\ep

\subsection{Viability and comparison property for multi-dimensional BSDEs}\label{App Viability}

We generalize in this paragraph some viability and comparison properties for multidimensional BSDEs in a closed convex cone $\Cc$ of $\R^{2m}$, whenever we add some reflections on the $Y$-component of the BSDE. The two following propositions are respectively extensions of Theorem 2.5 in \cite{buqura00}  and a simplifying version of Theorem 2.1 in \cite{hp06}. Their derivations do not present major difficulty and we choose to detail them for sake of completeness. \\

Let $(Y,Z)$ $\in$ $(\Sc^{2}_{\F}\times\mathbf{L^2_{\F}(W)})^{2m}$ satisfying
\beq\label{BSDEviab}
Y_{t} & = & Y_{T}+\int_{t}^TF(s,Y_{s},Z_{s})ds-\int_{t}^T\langle Z_{s},dW_{s}\rangle+K_{T}-K_{t} \;, \quad 0\le t \le T\;,
\enq
where $F~:\Omega\times[0,T]\times\R^{2m}\times\R^{2m\times d} \rightarrow \R^{2m}$ is a progressively measurable function satisfying \textbf{(H2)} (ii) and $K$ is an  $\R^{2m}$-valued finite variation process such that
\beqs
K_{t} & = & \int_{0}^tk_{s}d|K|_{s}\;,
\enqs
with $k_{t}$ $\in$ $\Cc$ and $|K|_{s}$  the variation of $K$ on $[0,s]$. We denote by $d_{\Cc}$ the distance to $\Cc$, i.e. $d_{\Cc}:x\mapsto\min_{y\in \Cc}|x-y|$, and introduce $\Pi_{\Cc}$ the projection operator onto $\Cc$.

\begin{Proposition}\label{BVP}
Suppose $Y_{T}$ $\in$ $\Cc$ and there exists a constant $C^0$ such that $F$ satisfies
\beq\label{ApendixCondition}
4\langle y-\Pi_{\Cc}(y),F(t,y,z)\rangle & \leq & \langle D^2|d_{\Cc}|^2(y)z, z\rangle + 2C^0 |d_{\Cc}|^2(y)\, \qquad \P-a.s.\,,
\enq
for any $(t,y,z)\in[0,T]\times\R^{2m}\times\R^{2m\times d}$ such that  $|d_{\Cc}|^2$ is twice differentiable at the point $y$.
Then, we have 
\beqs
Y_{t} & \in & \Cc\,, \qquad 0 \le t\le T \;, \qquad \P-a.s.
\enqs
\end{Proposition}

\ni\textbf{Proof.} 
The proof presented here is an adaptation of the one of Theorem 2.5 in \cite{buqura00}, allowing to tackle the additional difficulty due to the $dK$ term in the dynamics of $Y$. \\ Let $\eta$ $\in$ $C^\infty(\R^{2m})$ be a non-negative function, with support in the unit ball, such that $\int_{\R^{2m}}\eta(x)dx$ $=$ $1$.
For $\delta$ $>$ $0$ and $x$ $\in$ $\R^{2m}$, we define 
\beqs
\eta_{\delta}(x)  := \frac{1}{\delta^{2m}}\eta\Big(\frac{x}{\delta}\Big) &\quad\mbox{ and }\quad&
\phi_{\delta} (x) := 
\int_{\R^{2m}}|d_{\Cc}(x-y)|^2\eta_{\delta}(y)dy \;.
\enqs 
Via direct computation, one can verify that 
$\phi_{\delta}$ $\in$ $C^\infty(\R^{2m})$ and 
\begin{equation}\label{PropPhi}\left\{\begin{array}{l}
0 ~\leq~ \phi_{\delta}(x) ~\leq~ (d_{\Cc}(x)+\delta)^2 \;,
\\
D\phi_{\delta}(x)~=~\int_{\R^{2m}}D|d_{\Cc}(y)|^2\eta_{\delta}(x-y)dy ~~ \mbox{and} ~  |D\phi_{\delta}(x)|~\leq~2(d_{\Cc}(x)+\delta)\;,\\
D^2\phi_{\delta}(x)~=~\int_{\R^{2m}}D^2|d_{\Cc}(y)|^2\eta_{\delta}(x-y)dy ~ \mbox{and} ~~ 0~\leq~  |D^2\phi_{\delta}(x)|~\leq~2I_{2m} \;,
\end{array}\right.
\end{equation}
for any $x\in\R^{2m}$. An application of It\^{o}'s formula to $\phi_{\delta}(Y)$, combined with these estimates and $d_{\Cc}(Y_T)=0$, leads to 
\beq\label{controlpropphi}
\Eb\phi_{\delta}(Y_{t}) & = &  \Eb\phi_{\delta}(Y_{T}) + \Eb \int_{t}^T \langle D\phi_{\delta}(Y_{s}), F(s,Y_{s},Z_{s})\rangle ds-\frac{1}{2} \Eb\int_{t}^T\langle D^2 \phi_{\delta}(Y_{s})Z_{s},Z_{s} \rangle ds \nonumber \\ 
& & + \Eb\int_{t}^T\langle D\phi_{\delta}(Y_{s}), k_{s}\rangle d|K|_{s} \nonumber \\
 & \leq & \delta^2+ \Eb \int_{t}^T\hspace{-2mm} \int_{\R^{2m}}\Big[\langle D|d_{\Cc}(y)|^2, F(s,y,Z_{s})\rangle -\frac{1}{2} \langle D^2|d_{\Cc}(y)|^2Z_{s},Z_{s} \rangle\Big]\eta_{\delta}(Y_{s}-y)dy ds \nonumber\\ 
  & & -\Eb\int_{t}^T\hspace{-2mm} \int_{\R^{2m}}\langle D|d_{\Cc}(y)|^2, F(s,y,Z_{s})-F(s,Y_{s},Z_{s})\rangle\eta_{\delta}(Y_{s}-y)dy ds \nonumber\\
   & & +\Eb\int_{t}^T\hspace{-2mm} \int_{\R^{2m}}\langle D|d_{\Cc}(y)|^2, k_{s}\rangle\eta_{\delta}(Y_{s}-y)dy d|K|_{s}\,,  \qquad 0\le t\le T\;.
\enq
Since $k$ is valued in the closed convex cone $\Cc$, we observe that 
\beqs 
\langle D|d_{\Cc}(y)|^2, k_{s}\rangle &\leq& 0\;, \qquad 0\le s\le T\;, \quad y\in\R^{2m}\;.
\enqs 
Then, plugging this expression, \reff{ApendixCondition} and inequality $2d_c(.)\leq 1 + d_c(.)^2$ in \reff{controlpropphi},  we get
\beqs
\Eb\phi_{\delta}(Y_{t}) \!\!\!\!& \leq &\!\! \!\! \delta^2 + C^0 \Eb \int_{t}^T\hspace{-2mm}\int_{\R^{2m}}  |d_{\Cc}(y)|^2\eta_{\delta}(y-Y_{s}) dyds\\
 & & +2 \Eb\int_{t}^T\hspace{-2mm}\int_{\R^{2m}}d_{\Cc}(y)\eta_{\delta}(Y_{s}-y)  \max_{y':\;|y'-Y_{s}|\leq\delta}|F(s,y',Z_{s})-F(s,Y_{s},Z_{s})|dyds\\ 
 \!\! \!\!& \leq &\!\!\!\! \delta^2+C^0\!\int_{t}^T\!\!\Eb \phi_{\delta}(Y_{s})ds + \Eb\!\int_{t}^T\!\!\!(1+\phi_{\delta}(Y_{s}))\!\!\max_{y':\;|y'-Y_{s}|\leq\delta}|F(s,y',Z_{s})-F(s,Y_{s},Z_{s})|ds\,,
\enqs
 for any $t\in[0,T]$. Using the uniform Lipschitz property of $F$, we deduce
\beqs
\Eb \phi_{\delta}(Y_{t}) & \leq & C\left\{\delta^2 + \delta+ \int_{t}^T\Eb \phi_{\delta}(Y_{s})ds \right\} \;,\quad 0\leq t\leq T\;, \;\, \delta>0\,,
\enqs
and Gronwall's lemma leads to 
\beqs
\Eb\phi_{\delta}(Y_{t}) & \leq & C(\delta^2+\delta),\quad 0\leq t\leq T\;,\;\, \delta>0\,.
\enqs
Finally, from Fatou's Lemma, we have 
\beqs
\Eb |d_{\Cc}(Y_{t})|^2 & \leq &  \liminf_{\delta\rightarrow0}\Eb\phi_{\delta}(Y_{t}) ~=~0\;,\quad 0\leq t\leq T\;,
\enqs
which concludes the proof.
\ep

\vspace{2mm}

We now turn to the obtention of a multidimensional comparison result for BSDEs, whenever the dominating BSDE suffers additional reflections. This proposition also simplifies the results of Theorem 2.1 in \cite{hp06} in the case where the $i^{th}$ component of each driver only depends on the $i^{th}$ component of $Z$, for any $i\le d$. \\

Consider $(Y^1,Z^1,K^1)$ $\in$ $(\Sc^{2}_{\F}\times\mathbf{L^2_{\F}(W)}\times\mathbf{A^{2}_{\F}})^m$ satisfying
\beqs
Y_{t}^1 & = & Y^1_{T}+\int_{t}^TF_{1}(s,Y_{s}^1,Z^1_{s}) ds-\int_{t}^T\langle Z^1_{s},dW_{s}\rangle+K^1_{T}-K^1_{t}\,, \quad 0\le t\le T\;,
\enqs
and  $(Y^2,Z^2)$ $\in$ $(\Sc^{2}_{\F}\times\mathbf{L^2_{\F}(W)})^m$ satisfying
\beqs
Y_{t}^2 & = & Y^2_{T}+\int_{t}^TF_{2}(s,Y_{s}^2,Z^2_{s}) ds-\int_{t}^T\langle Z^2_{s},dW_{s}\rangle \;,
\quad 0\le t\le T\;,
\enqs
 where $F_1$ and $F_2$ are two driver functions satisfying {\bf (H2)} (ii) and such that the $i^{th}$ component of each driver only depends on the $i^{th}$ component of the corresponding $Z$, for any $i\le d$. 

\begin{Proposition}\label{thcomAp}
Suppose $Y^1_{T}$ $\geq$ $Y^2_{T}$ and the existence of a constant $C^1$ such that 
\beq\label{BVP eq}
- 2 \langle y,F_{1}(t,y',z)-F_{2}(t,y',z')\rangle & \leq & C^1|y|^2 + \sum_{i=1}^m |z_{i}-z_{i}'|^2~\mathbf{P}-a.s. \;, \quad\qquad
\enq
for any $(t,y,y',z,z')$ $\in$ $[0,T]\times(R^+)^m\times\R^m\times[\R^{m\times d}]^2$. Then $Y^1_{t}$ $\geq$ $Y^2_{t}$, for all $t$ $\in[0,T]$. 
\end{Proposition}

\textbf{Proof.}
The process $(Y^1-Y^2,Y^2)$ is valued in $\R^{2m}$ and solution of a BSDE of the form \reff{BSDEviab} associated to the driver
\beqs
F \; : \; (t,(y,y'),(z,z')) &\mapsto& ( F_1(t,y+y',z+z') -  F_2(t,y',z') ,  F_2(t,y',z'))\;,
\enqs
for any $\{t,(y,y'),(z,z')\}\in[0,T]\times\R^{2m}\times\R^{2m\times d}$.  Introducing the closed convex cone $\Cc$ $:=$ $(\R^+)^m\times\R^m$ of $\R^{2m}$, we see that 
$d_{\Cc}(y,y')=|y^-|$ for $(y,y')\in\R^{2m}$. Therefore, we deduce from the Lipschitz property of $F_1$ and \reff{BVP eq} that 
\beqs
  &&4\langle (y,y') - \Pi_{\Cc}(y,y'),F(t,(y,y'),(z,z'))\rangle \\
&=& 4\langle - y^-,F_{1}(t,y+y',z+z')-F_{1}(t,y',z+z')\rangle  + 4\langle - y^-,F_{1}(t,y',z+z')-F_{2}(t,y',z')\rangle \\
&\leq&  4 k |y^-|^2  +  2\sum_{i=1}^m\mathbf{1}_{y_{i}<0}|z_{i}|^2+2C^1|y^-|^2  \\
&= & \langle D^2|d_{\Cc}|^2(y,y')(z,z'), (z,z')\rangle + (2C^1+4k)  |d_{\Cc}|^2(y,y')\, \;\quad \P-a.s.\,,
\enqs
for any $\{t,(y,y'),(z,z')\}\in[0,T]\times\R^{2m}\times\R^{2m \times d}$. Applying Proposition  \ref{BVP} with $C^0=C^1+2k$, we deduce that the process $(Y^1-Y^2,Y^2)$ is valued in $\Cc$ and complete the proof.
\ep

\end{document}